\documentclass{article}
\usepackage{amsfonts,amsmath, amssymb}

\usepackage{xcolor}

\usepackage{enumerate}
\usepackage{epsf,epsfig,amsfonts,a4wide}
\usepackage{amsfonts}
\usepackage{amsmath,amssymb,amsthm}
\usepackage{url}
\usepackage{amsmath,amssymb,amsthm}
\usepackage{amssymb}
\usepackage{amsthm}
\usepackage{epsf,epsfig,amsfonts,graphicx}
\usepackage{a4wide}
\numberwithin{equation}{section}

\newtheorem{problem}{Problem}[section]
\newtheorem{question}[problem]{Question}
\newtheorem{con}[problem]{Conjecture}
\newtheorem{Th}[problem]{Theorem}

\newtheorem{proposition}[problem]{Proposition}

\theoremstyle{remark}

\newcommand{\be}{\begin{equation}}
\newcommand{\ee}{\end{equation}}

\newcommand{\tr}{\mathrm{tr}\,}

\newcommand{\R}{\mathbb{R}}\newcommand{\Id}{\textrm{\rm Id}}

\newcommand{\goth}{\mathfrak}

  \newcommand{\bbR}{\mathbb R}

  \newcommand{\bbC}{\mathbb C}

   \DeclareMathOperator {\Ad}    {Ad}

   \newcommand{\Fg}{\mathfrak g}

 \newcommand{\cA}{\mathcal A}

\newcommand{\NO}{no.~}

  \DeclareMathOperator {\rank}  {rank}
 
  \DeclareMathOperator {\corank}  {corank}

\newcommand{\weg}[1]{}

\title{Open Problems, Questions, and Challenges  in Finite-Dimensional Integrable Systems. }
\author{Alexey Bolsinov\footnote{
Department of Mathematical Sciences,
 Loughborough University,
 LE11 3TU, UK  and  Faculty of Mechanics and Mathematics, Moscow State University, Moscow 119991, Russia \ \
 \quad {\tt A.Bolsinov@lboro.ac.uk}},   Vladimir S. Matveev\footnote{Institut f\"ur Mathematik,
Fakult\"at f\"ur Mathematik und Informatik,
Friedrich-Schiller-Universit\"at Jena,
07737 Jena, Germany \ \  \quad  {\tt vladimir.matveev@uni-jena.de}}, Eva Miranda\footnote{{Department of Mathematics}, Universitat Polit\`{e}cnica de Catalunya and BGSMath, Barcelona, Spain, {\tt eva.miranda@upc.edu}}, Serge Tabachnikov\footnote{Department of Mathematics,
Pennsylvania State University, University Park, PA 16802, USA   {\tt tabachni@math.psu.edu}}
 }
 \date{}

\begin{document}
\maketitle

\begin{abstract}
The paper surveys open problems and questions related to different aspects of integrable systems with finitely many degrees of freedom. Many of the open problems were  suggested by the participants of the conference ``Finite-dimensional Integrable Systems, FDIS 2017'' held at CRM, Barcelona in July 2017.
\end{abstract}


\tableofcontents

\section{Introduction}

Many different and formally nonequivalent definitions of finite dimensional integrable systems will be used in this paper.
The basic   one  is as follow: we consider a symplectic manifold $(N^{2n}, \omega)$, a function $H:N^{2n}\to \mathbb{R}$ on it,  and the  Hamiltonian system corresponding to it.   It is {\it integrable}, if there exist
$n-1$ function $F_1,...,F_{n-1}$ (called integrals) such that the  $n$-functions $H, F_1,...,F_{n-1}$ commute with respect to the Poisson structure corresponding to the symplectic structure, and their differentials are linearly independent at almost every point of the manifold.  Different generalization of this basic definition will be of course considered, for example, for  Poisson manifolds (for a Poisson structure of  rank $2n$ on a $2n +r$ dimensional manifold we require the existence of $r$ Casimirs  $F_{n},...,F_{n+r-1}$  and $n$  commuting integrals, see \S \ref{bi} for details).    We will also consider the case  when the Hamiltonian system is discrete, i.e., just a symplectomorphism of $N^{2n}$; then we require the existence of $n$ functions which commute  with respect to the Poisson structure corresponding to the symplectic structure and which are preserved by the symplectomorphism.

The theory of integrable systems  appeared as a group of methods  which can be applied to find  exact solutions of dynamical systems, and has grown in a separate field of investigation, which uses methods from other branches of mathematics and has applications in different fields.

The present survey collects open problems in the field of finite dimensional integrable systems, and also   in
other fields, whose solution may influence the theory of integrable systems.
We do    hope that such a list will be useful in particular  for early stage researchers, since many problems in the list can, in fact, be seen as  possible topics of  PhD theses,
  and also  for people from other branches of mathematics,    since for many problems in our list their  solutions are  expected to use `external' methods.

  Many open problems in our list are fairly complicated and are possibly out of immediate  reach; in this case we tried to also present easier special cases. But even  in these cases we tried to keep the level of problems high -- in our opinion, for most problems in our list, the
	solution deserves to be published.

This list is motivated by the series of FDIS conferences organized every second year at different locations (in 2017, in Barcelona),  and which is possibly  the most important  series of conferences in the theory of finite dimensional integrable systems.  Many problems in our  list are suggested by the participants of this  conference. Needless to say, the final selection of the material reflects the taste, knowledge, and interests of the authors of this article.

\paragraph{Acknowledgements.} We are grateful to our colleagues who responded to our request to share open problems and conjectures on finite-dimensional integrable systems with us: M. Bialy, D. Bouloc,   H. Dullin,  A. Glutsyuk, A. Izosimov, V. Kaloshin, R. Kenyon, B. Khesin, A. Pelayo, V. Retakh,  V. Roubtsov,  D. Treschev, S. V\~u Ng{\d o}c, and N. Zung. We thank the referees of this paper for their useful comments and suggestions.

A. B. is supported by the Russian Science Foundation (project no. 17-11-01303).
  V. M. is supported by the DFG. E. M. is supported by the Catalan Institution for Research and Advanced Studies via an ICREA Academia Prize 2016 and partially supported  by the grants MTM2015-69135-P/FEDER and  2017SGR932 (AGAUR).
	S.T. was supported by NSF grant  DMS-1510055.

\section{Topology of integrable systems,  Lagrangian fibrations, and their invariants}

Recall that a finite dimensional   integrable system on a symplectic manifold $(M,\omega)$ gives rise to a {\it singular
Lagrangian fibration} whose fibres, by definition, are common level manifolds of the integrals $F_1,\dots, F_n$. Our assumption that
 the integrals commute implies that nondegenerate   fibers are  Lagrangian submanifolds, i.e., their dimension is half the dimension of the manifold, and the restriction of $\omega$ to them  vanishes. The word `singular' means that we allow fibers which are not manifolds. Such fibers typically contain  points $x\in M$ where the differentials $dF_1(x), \dots, dF_n(x)$ are linearly dependent.

Consider a singular Lagrangian fibration $\phi : M \to B$  with compact fibers  associated with a finite-dimensional integrable system. In general, the base $B$ of this fibration is not a manifold, but can be understood as a stratified manifold whose strata correspond to different ``singularity levels''.  If we remove all singular fibers from $M$, then the remaining set $M_{\mathrm{reg}}$ is foliated into regular Liouville tori, so that we have a locally trivial fibration  $M_{\mathrm{reg}} \overset{T}{\longrightarrow} B_{\mathrm{reg}}$ over the regular part $B_{\mathrm{reg}}$ of  $B$. Notice that  $B_{\mathrm{reg}}$ is the union of strata of maximal dimension and, in general, may consist of several connected components.

Recall that in a neighbourhood of each regular fiber (Liouville torus), one can define action-angle variables.  The actions can be treated as local coordinates on $B$.  Since they are defined up to integer affine transformations,  they  induce an integer affine structure on $B_{\mathrm{reg}}$.   In a neighbourhood of a singular stratum,  this affine structure, as a rule, is not defined but we still may study its asymptotic behaviour and think of it as a singular affine structure on $B$ as a whole.

The action variables are preserved under symplectomorphisms and, therefore, we may think of them as symplectic invariants of  singular Lagrangian fibrations.  How much information do they contain?  Examples show that very often the action variables (equivalently, singular affine structure on $B$) are sufficient to reconstruct completely the structure of the Lagrangian fibration up to symplectomorphisms.  The classical result illustrating this principle is Delzant theorem \cite{Delzant}  stating that, in the case of Lagrangian fibrations related to Hamiltonian torus actions, the base $B$ is an affine polytope (with some special properties) which determines the structure of the Lagrangian fibration (as well as the torus action) up to a symplectomorphism.  In particular, two Lagrangian fibrations $M^{2n} \overset{T^n}{\longrightarrow} B$ and $\tilde M^{2n} \overset{T^n}{\longrightarrow} \tilde B$ are symplectically equivalent if and only if their bases $B$ and $\tilde B$ are  equivalent as spaces with integer affine structures.

A similar statement still holds true for some Lagrangian fibrations with more complicated singularities  (both non-degenerate and degenerate), see \cite{Toulet,  San, Parabolic}.  We believe that an analog of Delzant theorem holds true in a much more general situation.  Ideally, such an analog could be formulated as the following principle  (it is not a theorem, as a counterexample is easy to construct!):

{\it Let $\phi: M \to B$ and $\phi':M' \to B'$ be two singular Lagrangian fibrations.  If $B$ and $B'$ are affinely equivalent (as stratified manifolds with singular affine structures), then these Lagrangian fibrations are fiberwise symplectomorphic.  }

\begin{question}
Under which additional assumptions  does this principle become a rigorous theorem?
\end{question}

In our opinion, this is an important general question, which is apparently quite difficult to answer in full generality.   Some more specific questions could be of interest too.

Let $\phi: M^4 \to B$ be an almost toric fibration (see \cite{Symington2, Symington1}),  which means, in particular, that its singularities are all non-degenerate and can be of elliptic and focus type only.  Consider a  typical  situation when the base $B$ of such a fibration is a two-dimensional domain with boundary (having some `corners')  and some isolated singular points    of the focus type.   This domain is endowed with an integer affine structure, having singularities at focus points.

\begin{problem}
Consider two Lagrangian fibrations $\phi: M^4 \to B$ and $\phi': {M'}^4 \to B'$.   Assume that $B$ and $B'$ are affinely equivalent in the sense that there exists an affine diffeomorphism $\psi : B \to B'$.  Is it true that under these assumptions the corresponding Lagrangian fibrations are symplectomorphic?
\end{problem}

In the local and semi-global setting, we may ask a similar question for a Lagrangian fibration in a neighbourhood of a singular point or of a singular fiber.  According to  the   Eliasson theorem  \cite{Eliasson, mirandathesis, mirandazung}, in the non-degenerate case, there are no local symplectic invariants.  In other words,  non-degenerate singular points of the same topological (algebraic) type are locally fiberwise symplectomorphic.  The case of non-degenerate singular fibers has been discussed in \cite{Toulet, San, Pelayo2, Pelayo3} for basic singularities of elliptic, hyperbolic, and focus types.  One of the main observations made in these papers is that every action variable can be written in the form $I = I_{\mathrm{sing}} + I_{\mathrm{reg}}$, where the singular part $I_{\mathrm{sing}}$ is the same for all singularities  of a given topological type (due to Eliasson's result) and  $I_{\mathrm{reg}}$ is a smooth function which can be arbitrary.   This regular part (or its Taylor expansion) in many cases may serve as a complete symplectic invariant in semi-global setting.

Let us now discuss what happens in the case of  degenerate   singularities.

\begin{problem}
Do local symplectic invariants exist for diffeomorphic degenerate singularities? How many and of what kind are they?  This question makes sense even in the simplest case of one-degree-of-freedom systems.
\end{problem}

As we know from concrete examples of integrable systems, many degenerate singularities are stable in the sense that they cannot be destroyed by a small perturbation of the system \cite{BF, Giacobbe, Kalash}.  The problem of description of stable degenerate singularities is very important on its own right, but for many degrees of freedom,  it is quite complicated.  However,  in the case of two degrees of freedom  (perhaps depending on a parameter), there is a number of well-known examples of stable bifurcations (e.g.,  integrable period-doubling \cite{BF}  and Hamiltonian Hopf bifurcation \cite{dp}).  It would be interesting to clarify the situation with asymptotic behaviour of angle variables at least for them.

 V.~Kalashnikov  \cite{Kalash} has  described all stable singularities of rank 1 for Hamiltonian systems of two degrees of freedom.

\begin{problem}
Describe the symplectic invariants of stable singularities from \cite{Kalash}.  For such singularities, one of  the action variables, say $I_1$, is smooth, and the other $I_2$ is singular. Is it true that  $I_1$ and $I_2$ are sufficient for the symplectic classification?

Here is an almost equivalent version of this question. Assume that we have two functions $H$ and $F$ that commute with respect to two different symplectic structures $\omega_1$ and $\omega_2$  and define (in  both cases) a stable singularity from \cite{Kalash}.   Assume that the action variables $I_1$ and $I_2$ are the same for $\omega_1$ and $\omega_2$. This condition is equivalent to the relation
$$
\oint_\gamma \alpha_1 = \oint_\gamma \alpha_2
$$
for any  cycle $\gamma$ on any regular fiber $\mathcal L_{f,h}=\{ F=f, H=h\}$ and appropriately chosen $1$-forms $\alpha_i$ such that $d\alpha_i = \omega_i$, $i=1,2$.  Is it true that, under these conditions, there exists a smooth map $\psi$ that preserves the functions $F$ and $H$ and such that
$\psi^*(\omega_2)=\omega_1$? In the simplest case the answer is in \cite{Parabolic}.

\end{problem}

\begin{question} \label{prob:1.4}
 Assume that we know explicit formulas for the action variables $I_1,\dots, I_n$ so that we are able to analyse their asymptotic behaviour  in a neighbourhood of a singular fiber.   Can we recover the topology of this singularity from the asymptotic behaviour  (or at least to distinguish between different types of singularities)?  For example, we know that, in the case of non-degenerate hyperbolic singularities,  the singular part $I_{\mathrm{sing}}$ is of the form $h\ln h + \dots$ . Is this property characteristic for non-degenerate  hyperbolic singularities? The case of one degree of freedom is understood in \cite{phdlorenzo}.
\end{question}

For degenerate singularities, this question becomes, of course, more interesting and important.

Another interesting problem  would be a description of singularities, both in local and semi-global setting, which may occur in algebraically integrable systems. There is no commonly accepted definition of algebraic integrability  (see discussion on this subject and relevant definitions in \cite{Vanhaecke_LN}) and the answer may depend on the approach used. In the most natural setting,  algebraic integrability assumes that the commuting functions $F_1,\dots, F_n$ are polynomials and one can naturally complexify all the related objects  (including the phase space, symplectic or Poisson structure, and time) so that the integrable system under consideration becomes the real part of a certain complex integrable system. The crucial additional condition on this complex system is that its general complex fibers are isomorphic to affine parts of Abelian varieties (see \cite{Vanhaecke_LN}) which imposes rather strong topological conditions on singularities.
For instance, in the case of one degree of freedom, it implies that the Milnor number of a singular point cannot be greater than 2.  In other words, singularities of algebraically integrable systems are expected to be very special. 

\begin{problem} Describe all the topological types of singularities that may appear in algebraically integrable systems with a small  {\rm(}$\le  3${\rm)} degree of freedom. Next,  describe the symplectic invariants of such singularities.
\end{problem}

We expect that  solving local and semi-global problems mentioned above will simplify the passage to global results.  Some  open problems related to classification of integrable systems in global setting are stated below.

Recall that an integrable system $F=(J,H)$ is called {\it semitoric} if all of its singularities are nondegenerate and of either elliptic or focus types
 (hyperbolic singularties are forbidden) and, in addition, $J$ is proper and  generates a Hamiltonian $S^1$-action on the manifold.

\begin{problem}[\'{A}. Pelayo]
Extend the classification of semitoric systems $F=(J,H)$ in \cite{Pelayo2, Pelayo3} to allow for $F$ having   non\--degenerate singularities with
hyperbolic blocks  (a version of this problem was mentioned in various forms in~\cite{DCDS12}).
\end{problem}

The first step in this direction is in Dullin {\it et al} \cite{dp}, see also a general approach to the symplectic classification problem in Zung  \cite{Zung}.

Note that,  in contrast with the semitoric case, even the topology of a  small  neighbourhood  of a singular fiber can be quite complicated, see e.g. \cite{bolsinovmatveev, BolsinovIzosimov, IzosimovVestnik, lerman,  matveev1, matveev2, matveevoshemkov,  PelayoTang, Zung1,Zung0,Zung}.

There are also a number of interesting problems at the intersection of integrable systems and Hamiltonian compact group actions, the following being one of them.

\begin{problem}[\'{A}. Pelayo]
Consider a compact connected $2n$\--dimensional symplectic manifold $M$, endowed with a Hamiltonian $(S^1)^{n-1}$\--action; these are called
\emph{complexity one spaces}. Consider an integrable system $f_1,\ldots, f_{n}$ on $M$ where $(f_1,\ldots, f_{n-1}) \colon M \to \mathbb{R}^{n-1}$
is the momentum map  of the   Hamiltonian $(S^1)^{n-1}$\--action. Suppose that the singularities of the integrable system are non\--degenerate and
also they do not contain hyperbolic blocks. Study how the invariants of the complexity one space are related to the invariants
of the semitoric system.
\end{problem}

For  recent progress in  the case  $2n=4$ see   Hohloch {\it et al}  \cite{HoSaSe}.

Recently there have been many works on integrable systems which can be obtained
by the method of toric degenerations.
In particular, many well-known integrable systems can be constructed this way, including the bending flow of 3D polygons Foth {\it et al}  \cite{Foth-Hu2005}
and the Gelfand-Cetlin\footnote{Another common spellings are  Zetlin and Zeitlin.} system Nishinou {\it et al} \cite{Nishinou-Nohara-Ueda}.   See also \cite{HK2012, Kaveh, NU2014}.
 A common feature of
integrable systems coming from toric degenerations is that their base space (the space
of fibers of the momentum map) is a convex polytope, which implies that the singular fibers
are quite special.

\begin{problem}[N. T.  Zung] Study the topology and geometry of these singular fibers
and their small neighbourhoods.
\end{problem}

First results in this directions are due to
D.~Bouloc  \cite{BoulocSympGeom2016}, where
the bending flows of 3D polygons and the  system on the 2-Grasmannian manifold were studied, and by Bouloc {\it et al}  \cite{boulocmirandazung} for the
 Gelfand-Cetlin system.

The natural  conjecture is that, in most  cases, the degenerate singular fibers are still
smooth manifolds or orbifolds and, moreover, they are isotropic (or even Lagrangian).  They represent a new interesting class of singularities in integrable systems for which almost all traditional questions remain open (symplectic invariants, normal forms, stability).

Another problem is related to the so called SYZ model of mirror symmetry.
In the case of Lagrangian fibrations on Calabi-Yau manifolds (see for example \cite{matessi, GrossSiebert1, GrossSiebert2}) of complex dimension 3 (real dimension 6)
there are 1-dimensional families of singularities of focus-focus type (which are well-understood
from the point of view of integrable systems), but there is also, for each system, a finite number
of special singular fibers ``at the intersection'' of these families of focus-focus singularities, and these
special singular fibers are still a mystery.

\begin{problem}[N. T. Zung] Give a clear description of these
special singular fibers.
\end{problem}


\section{Natural Hamiltonian systems on closed manifolds}


By a natural  Hamiltonian system we understand a Hamiltonian system constructed as follows. Take     a
manifold $M=M^n$ equipped with a metric $g=g_{ij}$ and a function $U:M\to \mathbb{R}$, then consider the  Hamiltonian system   on the cotangent bundle $T^*M$  with the Hamiltonian given  by
\begin{equation} \label{ham}
H= K + U = \tfrac{1}{2} g^{ij} p_i p_j + U.
 \end{equation}
 Here $p=(p_1,\dots, p_n)$ denote the momenta (i.e., the coordinates on $T^*_xM$ associated to local  coordinates $x=(x_1,...,x_n)$ on $M$) and $g^{ij} $ is the inverse matrix to $g_{ij}$. Summation over repeating indices is assumed, and the symplectic structure on $T^*M$ is canonical.

Of course, $H$ itself is an integral for such a system, which we call trivial. Also all functions of  $H$  will be considered as trivial integrals.

\subsection{  Integrability in the class of integrals  that  are polynomial in momenta. }

 In this section we  discuss the existence of an integral $F$ (or a family of integrals $F_1,...,F_{n-1}$)   which is polynomial of some fixed degree in momenta.  The coefficients of this polynomial can, and usually do, depend on the position of the point $x\in M$, and in local coordinates they are  functions on $M$.  By \cite{kruglikov2016}, if the manifold $M$ and the Hamiltonian $H$ are smooth,  these functions are automatically smooth too.

 This is one of the oldest  part of the theory of integrable systems. Indeed,  most, if not all,   classically known   integrable systems admit integrals which are polynomial in momenta.  For brevity, in what follows we will sometimes refer to them as {\it polynomial integrals}, meaning {\it polynomial in momenta}.

A generic natural integrable system  admits no nontrivial polynomial integral even locally, see \cite{kruglikov2016}. The existence of such an integral is therefore a nontrivial  differential-geometrical condition  on the metric $g$ and  the function $U$, which actually appears naturally  in some topics in  differential geometry. The relation to differential geometry becomes even clearer in the important special case,  when the function $U$ is identically zero.
   Such Hamiltonian systems are called  \emph{geodesic flows}.  It is well known, and is easy to check, that polynomial  integrals of geodesic flows are in natural one-to-one correspondence with the so-called Killing tensors.

    Note however that the condition  that $U$ is zero does not make the investigation of integrability easier. Indeed,
       by Maupertuis's principle (see e.g. \cite{Maup}), for each
        fixed value $h_0$  of the   energy $H=K+U$, the restriction of the Hamiltonian system to the isoenergy level $\{(x,p)\in T^*M \mid H(x,p)= h_0\}$ is closely related (= is the same up to a re-parameterisation) to  the geodesic flow of the metric  $g_{\mathrm{new}}=(h_0- U(x)) g$. In the latter formula we assume that $h_0\ne U(x)$, which can be achieved, for example, if the energy value $h_0$ is big  enough. Polynomial in momenta integrals for $H$  generate  then  polynomial in momenta integrals for the geodesic flow of $g_{\mathrm{new}}$.

 \subsubsection{Two-dimensional case.} \label{2}

By   Kolokoltsov \cite{Kolokoltsov1983}, no metric on a closed surface of genus $\ge 2$ admits any nontrivial polynomial integrals.  This result implies that no  natural Hamiltonian system on a  closed surface of genus $\ge 2$ admits such integrals either.
 If the surface is  the sphere  or  the torus,  the state of the art  is presented in   the following table:

\begin{center}
\begin{tabular} {|c|c|c|}
\hline
         & Sphere $S^2$  &Torus $T^2$ \\
         \hline
Degree 1 & All is known & All is known \\ \hline
 Degree 2 & All is
known & All is known\\ \hline
 Degree 3 & Series of examples & Partial negative results \\
\hline Degree 4 & Series of examples & Partial negative results \\
\hline Degree $\ge$ 5 & Almost  nothing is known  & Almost nothing is known \\
\hline
\end{tabular}

\end{center}

\vspace{3ex}

Here
``Degree''  means the smallest degree of a nontrivial  polynomial integral.
``All is known'' means that there exists an effective description and classification
(which can be derived, for example,  from  \cite{Bo-Ma-Fo1998}).
``Series of examples'' start with those coming from the Kovalevskaya and Goryachev-Chaplygin cases in rigid body dynamics, see e.g. \cite{Maup}, and  include  those constructed in \cite{dullin,kiyohara,valent}.

``Partial negative results'' have a different nature. Let us mention here the  series of papers by  M.~Bialy  \cite{Bialy1987},   N.~Denisova {\it et al.} \cite{De-2012} and   Agapov {\it et al.} \cite{Ag}. These papers study integrability, in the class of polynomial integrals, of natural Hamiltonian systems on $T^2$  under the condition that the metric $g$ is flat. It is proved that the existence of a nontrivial polynomial integral of degree $ m \le 4$ implies the existence of a nontrivial integral which is either linear or quadratic in momenta. An analogous result was obtained by A.~Mironov in  \cite{Mironov2010}   for $m=5$ under the additional assumption of real analyticity.

As for ``Degree $\ge 5$'' line of the table, let us first mention the results of Kiyohara \cite{kiyohara}  who, for each $d\ge 3$,
 constructed a family of metrics on $S^2$  whose geodesic flows  admit a polynomial integral of degree $d$, but do not admit any nontrivial  integrals of degree $1$ or $2$. Unfortunately, for the examples constructed in \cite{kiyohara},  it is not clear whether they might have nontrivial integral of ``intermediate''  degrees $d'$ (i.e., $2<d'< d$) or not.  The integral $F$ of degree $d$ constructed in \cite{kiyohara} is irreducible, but one may still imagine that there exists an additional polynomial integral of a smaller degree which is independent of both $H$ and $F$.
 The point is that  the  metrics constructed in \cite{kiyohara}   are {\it Zoll metrics}, i.e., all of their geodesics are closed. This means that the  geodesic flows are  superintegrable (see the discussion in
  \S \ref{zoll}),  which implies that  additional integrals exist; it is not clear though whether they are polynomial in momenta and, in fact, we conjecture that it is not the case.

  \begin{problem} \label{polynomial}  Complete the above table: {\rm (1)} construct new examples of natural Hamiltonian systems on closed two-dimensional surfaces
  admitting polynomial integrals, describe and, if possible, classify them; {\rm (2)} prove, if possible, nonexistence of such integrals (perhaps under some additional assumptions).
  \end{problem}

Of course,  this general  problem is very complicated, and  a solution of certain special cases  (some of them are listed below)  is  interesting. Even a construction of one new integrable example on a closed surface deserves to be published!

Let us list some conjectures/special cases/additional assumptions which were studied in  relation to  Problem  \ref{polynomial}.

\begin{con}[\cite{Maup}] \label{con:1}    {\it If the geodesic flow of
a Riemannian metric on the torus $T^2$   admits a nontrivial
integral   that is  polynomial of degree 3 in momenta,  then the metric admits  a  Killing vector field.  } \end{con}

An easier version of this  conjecture, which is also  interesting,  is as follows:

\begin{con}  \label{con:1b}  {\it If  a natural Hamiltonian system  with a nonconstant potential
  on $T^2$  admits an
integral   that is  polynomial of odd degree  in momenta,  then it admits a linear integral. If it admits a nontrivial integral of even degree, then it admits a nontrivial quadratic integral. } \end{con}

Resent results  in this conjecture include \cite{Corsi}.
Even the case when the integral has degrees 3 and 4 is not solved!  Yet an easier version is the following conjecture:

\begin{con}[Communicated by M. Bialy] \label{con:2}   {\it  If  a natural Hamiltonian system  with a nonconstant potential
  on a flat two-torus $T^2$ {\rm(}i.e., $g$ has zero curvature{\rm)} admits a polynomial
integral   of odd  degree,  then it admits a linear integral.  If it admits a nontrivial integral of even degree, then it admits a nontrivial quadratic integral.} \end{con}

As explained above, the  cases $d=3,4,5$ of this conjecture were proved   in  \cite{Bialy1987, Ag,  Mironov2010,De-2012}.  If the potential $U$ is a trigonometric polynomial, the problem was solved by N. Denisova and V.~Kozlov   in \cite{KozDen}.

Analogs of Conjectures \ref{con:1}, \ref{con:1b},  \ref{con:2}    are also interesting if the metric $g$ is semi-Riemannian  (i.e., an indefinite metric of signature $(1,1)$), which, of course, implies that our surface is a torus or  a Klein bottle.  Note that semi-Riemannian metrics  on $T^2$ whose geodesic flows admit integrals which are linear or quadratic in momenta were described in \cite{japan} and, combining this result with \cite{pucacco}, one obtains a description of all natural Hamiltonian systems on semi-Riemannian 2-tori with linear and quadratic integrals.

As seen from the discussion above, in the case of the sphere $S^2$,  almost nothing is known about integrals of degree $d\ge 5$, which suggests the following open problem:

 \begin{problem} \label{two-sphere}
Construct  a   natural Hamiltonian system  with a nonconstant potential
  on $S^2$  that  admits a nontrivial polynomial
integral  of degree $5$   and does not admit any nontrivial integrals of smaller degrees.
\end{problem}

Additional interesting generalisations  of the above questions appear    if  we consider a slightly more general class of Hamiltonian systems on $T^*M$, namely such that
$$
H:=K + L+ U,
$$
where $K$ and $U$ are as in \eqref{ham}, and $L$ is a linear in momenta function (sometimes referred to as magnetic terms); the solutions of this Hamiltonian system describe trajectories of charged particles in a magnetic field.  In this case,
 only  few explicit examples are known (e.g., those coming from Kovalevskaya, Clebsh, and Steklov integrable cases in  rigid body dynamics, and also from Veselova system and Chaplygin ball from  non-holonomic integrable systems), so any new example would be an interesting result.   Notice that even an explicit local description of such systems in the simplest case of quadratic integrals in dimension two remains an open problem  (we refer to \cite{SokMar} by V.~Marikhin   {\it et al} and  \cite{pucacco2} by G. Pucacco   {\it et al}  for the latest progress in the area).

 Note that in this case one can consider integrability on the whole cotangent bundle, as it was done, for example, in \cite{SokMar}, or on one energy level as,   for example, in \cite[\S3]{BBM}; both types of integrability are interesting.

 We also would like to draw  attention to two recent  papers containing  open problems on integrable natural Hamiltonian systems on  two-dimensional manifolds: Burns {\it et  al} \cite{burns}, see \S 10.2 there,  and Butler \cite{butler}, see \S\S 3.3--3.5 there.

\subsubsection{ Polynomial in momenta integrals in higher dimensions. }

In most known examples, the integrals of natural systems are either related to Hamiltonian reduction leading to integrable systems on Lie algebras (which may possess nontrivial polynomial integrals of arbitrary degree, see e.g. \cite{Jovanovich}), or the integrals are quadratic in momenta and  the quadratic terms are simultaneously  diagonalisable. There are also examples with polynomial integrals of higher degrees coming from physical models, such as the Toda lattice and Calogero--Moser systems.
Generally,   one simply needs more examples, which either are interesting to mathematical physics or ``live'' on closed manifolds.

  \begin{problem} \label{thigher}
Construct  new examples of  natural Hamiltonian systems on higher dimensional manifolds which are integrable in the class of integrals polynomial in momenta.
\end{problem}

Another  interesting direction of research might be relevant to  General Relativity.
 Recall that the most famous solutions of the Einstein equation, Schwarzschild and Kerr metrics, have integrable geodesic flow.

\begin{question}  Construct stationary axially symmetric 4-dimensional Einstein metrics admitting Killing tensors of higher order.
\end{question}

The condition ``stationary axially symmetric'' simply means the existence of two commuting vector fields, one of which is space-like and the other time-like.

Let us also touch upon two problems related to the homogeneous case.
 Recall that a {\it scalar invariant  of order} $k$ of a Riemannian manifold $(M, g)$  is a function on $M$, constructed from
covariant derivatives  of order $k$ or less of the curvature tensor  $R$ or, equivalently, an algebraic expression in the components of covariant derivatives of $R$ (of order $\le k)$ which transforms as a function under coordinate changes.  We say  that  $(M,g)$ is $k$-{\it homogeneous}, if all its scalar invariants up to order $k$ are constants.

It is known that,  for a sufficiently big $k$, any $k$-homogeneous manifold is automatically locally homogeneous, i.e., admits a locally transitive isometry group. As it follows from Singer  \cite{singer}, $k=n(n+1)/2+1$, where $n$ is the dimension of the manifold, is sufficiently big in this sense; later better estimations were found in \cite{vanhecke, Olmos}.

\begin{question}[Gilkey \cite{gilkey}]
 In the Riemannian case, is every $1$-homogeneous manifold   locally homogeneous? \end{question}

  The  question  was positively answered  for dimensions $3,4$.
The natural analogs  of  this conjecture for metrics of non-Riemannian signature are wrong. Notice that, in terms of polynomial integrals,  local homogeneity is equivalent to the existence of linear integrals of the geodesic flow that at each point span the tangent space.

 \begin{question}  In  a
symmetric  space, is every Killing tensor  a sum of symmetric products of Killing vectors?  Equivalently, is it true that the algebra of all polynomial integrals of the geodesic flow in a symmetric space is generated by linear integrals? \end{question}

Note that for general metrics the answer is clearly negative:  there exist metrics admitting nontrivial Killing tensors and no Killing vectors. For constant curvature spaces, the answer is positive  \cite{thomson}.

\subsection{ Superintegrable systems}  \label{zoll}

 Recall that a natural Hamiltonian system is {\it superintegrable}, if it admits more functionally independent integrals (of some special form)   than the dimension of the manifold   (in addition, one usually assumes that common level surfaces of these integrals are isotropic submanifolds in $T^*M$).  A system is called {\it maximally superintegrable}, if the number of functionally independent  integrals is   $2 n -1$,   which is the maximal number for an $ n$-dimensional manifold.

 An interesting and well-studied
    class of superintegrable  natural Hamiltonian systems consists of
     geodesic flows of Zoll metrics. Recall that a Riemannian metric (on a closed manifold) is called a {\it  Zoll} metric if all  its geodesics are closed, see e.g. \cite{Besse}. In this case, the Hamiltonian flow can be viewed as a symplectic  action   of $S^1$, the orbit space is an orbifold, and any function defined on this orbifold is an integral. Notice that locally, near almost every point, each  natural Hamiltonian system is superintegrable.  It is therefore standard to consider superintegrability in the class of
     integrals of a special form, typically, polynomial in momenta.  Then, superintegrability is a nontrivial condition on the system,  even locally.

   In dimension two, superintegrable systems with quadratic integrals were locally   essentially described in \cite{Darboux}.  Superintegrable systems with one linear and one cubic integral were described in \cite{seva}; the methods can be generalised for the case when one integral is linear and the second of arbitrary degree, see e.g. \cite{Val2}.

   \begin{problem}   {\it  Construct a  natural Hamiltonian system on the 2-sphere with a nonconstant potential which is superintegrable by integrals of degree  $\ge 3$ and admits no nontrivial integral of  degree one and two. } \end{problem}

  In higher dimensions, partial  answers were obtained  for conformally flat metrics only and are  well elaborated  only in dimension 3  for quadratic integrals,  see e.g. \cite{kress1,kress2}. The ultimate
    goal is, of course, to describe all superintegrable (in the class of polynomial in momenta integrals) natural Hamiltonian  systems. As the first nontrivial problem we suggest the following question:

 \begin{question}  Does there exist  a  non-conformally flat metric on the  sphere $S^n$, $n>2$, whose geodesic flow is
maximally superintegrable (in the class of integrals which are polynomial in momenta)?
\end{question}


\section{Integrable maps and discrete integrability}


 \subsection{Definition of integrability} \label{defint}

Iterations of a smooth map $f: M \to M$ yield a discrete-time dynamical system on a manifold $M$. 
More precisely, the evolution of a point $x\in M$ is given by repeated applications of the map $f$:
$$
x \mapsto f(x) \mapsto f^{\circ 2}(x) \mapsto \ldots f^{\circ k}(x) \mapsto \ldots
$$
Thus the time takes non-negative integral values. If the map is invertible, then the orbit of a point is also defined for negative values of $k$.

By Liouville integrability of such a system one means the following:
\begin{itemize}
\item The manifold $M^{2n}$ carries an $f$-invariant symplectic structure $\omega$;
\item There exist smooth functions (integrals) $I_1,\ldots,I_n$ on $M$, functionally independent almost everywhere, invariant under the map $f$, and Poisson commuting with respect to the Poisson structure induced by $\omega$.
\end{itemize}
 See \cite{Veselov91} for definitions and examples.

A modification of this definition is given by a Poisson manifold $M^{p+2q}$ whose Poisson structure has corank $p$. Then Liouville integrability of a map $f: M \to M$ means that the Poisson structure is $f$-invariant, and $f$ possesses almost everywhere functionally independent integrals $I_1,\ldots, I_{p+q}$, of which first $p$ are Casimir functions (the functions that Poisson commute with everything).

A  recent popular example is the pentagram map introduced in 1992 by R. Schwartz \cite{Schw}, a projectively invariant iteration on the moduli space of projective equivalence classes of polygons in the projective plane. This map is Liouville  integrable \cite{OST1,OST2}; its algebraic-geometric integrability was established in \cite{Soloviev}. The  integrability of the pentagram map is closely related with the theory of cluster algebras, see, e.g., \cite{FoMa,GSTV,Glick,GP,GoKe}.

\subsection{Billiards and billiard-like systems}

The billiard system describes the motion of a free particle in a domain with elastic reflection off the boundary. Let
$\Omega\subset  \mathbb{R}^2$  be a strictly convex domain bounded by a smooth curve $\gamma$. A point moves inside $\Omega$ with unit speed until it reached the boundary where it reflects according to the law of geometrical optics:   the
angle of reflection equals the angle of incidence, see, e.g., \cite{Tabachnikov95}.

An extension to higher dimensions and, more generally, to Riemannian manifolds with boundary is immediate: at the moment of reflection, the normal component of the velocity changes sign, whereas the tangential component remains the same.

The phase space of the billiard system consists of oriented lines (or oriented geodesics), endowed with its canonical symplectic structure, obtained by symplectic reduction from the symplectic structure of the cotangent bundle of the ambient space. These oriented lines are thought of as  rays of light;  the billiard map takes the incoming  ray to the reflected one.

In dimension two, an invariant curve of the billiard map can be thought as a 1-parameter family of oriented lines. Their envelope is called a billiard {\it caustic}: a ray tangent to a caustic remains tangent to it after the reflection.

A variation is the billiard system in the presence of a magnetic field. In the plane, or on a Riemannian surface, a magnetic field is given by a function $B(x)$, and the motion of a charged particle is described by the differential equation $\ddot x = B(x)J\dot x$, where $J$ is the rotation of the tangent plane by $\pi/2$. This equation implies that the speed of the particle remains constant (the Lorentz force is perpendicular to the direction of motion). In the Euclidean plane, if the magnetic field is constant, the trajectories are circles of a fixed radius, called the Larmor circles.

In general, a magnetic field on a Riemannian manifold $M$ is a closed differential 2-form $\beta$, and the magnetic flow is the Hamiltonian flow of the usual Hamiltonian $|p|^2/2$ on the cotangent bundle $T^*M$ with respect to the twisted symplectic structure $dp \wedge dq + \pi^*(\beta)$, where $\pi : T^*M \to M$ is the projection of the cotangent bundle on the base.

Another generalization are the Finsler billiards: the reflection law is determined by a variational principle saying that the trajectories extremize the Finsler length (which  is not necessarily symmetric), see \cite{GT}. Magnetic billiards are particular cases of Finsler billiards, see  \cite{Tabachnikov04} or \cite{BHTZ}.

Still another variation is the outer billiard, an area preserving transformation of the exterior of a closed convex plane curve. Like inner billiards, outer ones can be defined in higher dimensions: a curve in the plane is replaced by a closed smooth strictly convex hypersurface in linear symplectic space, and the outer billiard transformation becomes a symplectic mapping of its exterior.

Outer billiards can be defined in two-dimensional spherical and hyperbolic geometries as well. In the spherical case, the inner and outer billiards are conjugated by the spherical duality. See \cite{DT} for a survey.

Let us also mention two other non-conventional billiard systems: projective billiards \cite{Tabachnikov97} and symplectic billiards \cite{AT}.

\subsubsection{Around the Birkhoff conjecture}

It is a classical fact (see for instance \cite{Tabachnikov95} and references therein) that the billiard in an ellipse is integrable: the interior of the table is foliated by caustics, the confocal ellipses. Likewise, the billiard inside an ellipsoid is Liouville integrable, see, e.g., \cite{DR,Tabachnikov2005}.

The following is known as the Birkhoff conjecture, see \cite{Por}.

 \begin{con}  \label{bircon}
If a  neighbourhood of the boundary of the billiard table $\Omega \subset  \mathbb{R}^2$ is foliated by caustics, then the billiard curve $\gamma$ is an ellipse.
 \end{con}

We also formulate a generalized conjecture: {\it if an open subset of the  convex  billiard table is foliated by caustics, then this table is elliptic}.

A polynomial version of Birkhoff's conjecture concerns the case when the integral of the billiard map is polynomial in the velocity.
In the case of outer billiards, one modifies these conjectures in an obvious way (since outer billiards are affine-equivariant, there is no difference between ellipses and circles, and circles are integrable due to symmetry).

The state of the art of the (smooth) Birkhoff conjecture is as follows. The special case  of the conjecture in which  the caustics of the billiard trajectories foliate  $\Omega$, except for one point, was proved by M. Bialy \cite{Bialy1993}; later he extended this result to the elliptic and hyperbolic planes \cite{Bialy2013}. Bialy also proved a version of this result for magnetic billiards in constant magnetic field on surfaces of constant curvature \cite{Bialy2010}.

Recently,  substantial progress toward Birkhoff's conjecture was made by V. Kaloshin and his co-authors \cite{Avila,Kaloshin,Huang}: they proved the Birkhoff conjecture for small perturbations of ellipses. The general case remains open, and it is one of the foremost problems in this field.

\begin{problem}
Assume that the exterior of an outer billiard table is foliated by invariant curves. Prove that the table is an ellipse.
\end{problem}

The integral-geometric methods used by Bialy in his proofs fail due to non-compactness of the phase space.

The following is a generalization of Birkhoff's conjecture.

 \begin{con}  \label{Descon} {\rm (S. Tabachnikov)}
Let $\gamma\subset\mathbb{RP}^2$ be a smooth closed
strictly convex curve, and let $\mathcal F$ be a smooth foliation by convex closed curves in its exterior  neighbourhood, including  $\gamma$ as a leaf.
The intersection of every tangent line to $\gamma$  with the leaves of red {$\mathcal F$} define a local involution on this line,  the point of tangency being a fixed point.
Suppose that these involutions are projective transformations for all tangent lines  of $\gamma$. Then $\gamma$ is an ellipse, and the foliation consists of ellipses that form a pencil.
 \end{con}

The converse statement, that the intersections of a line with the conics from a pencil define a projective involution on this line, is a classical Desargues theorem.

A particular case of the above conjecture, when the involutions are isometries, yields the Birkhoff conjecture for outer billiards. In the generalized form, this conjecture implies, via projective duality, a version of Birkhoff conjecture for projective billiards and hence for the usual inner billiards.

The above formulated conjectures and problems have algebraic versions, and much progress has been made in this area recently. For inner billiards, the polynomial Birkhoff conjecture was proved in \cite{BM17,Glutsyuk17}, including the cases of non-zero constant curvature, and for outer billiards in the affine plane -- in \cite{GS}.

\begin{problem}
Prove the algebraic version of Birkhoff conjecture for outer billiards in a non-Euclidean surface of constant curvature.
\end{problem}

The next problem is a generalization of the previous one.

\begin{problem}
Prove Conjecture \ref{Descon} it the case that the foliation admits (i) a rational, (ii) an algebraic first integral.
\end{problem}

Algebraic non-integrability of magnetic billiards was studied by Bialy and Mironov in \cite{BM16}: they proved that a non-circular magnetic Birkhoff billiard is not algebraically integrable for all but finitely many values of the magnitude of the (constant) magnetic field. It is plausible that the result holds for every value of the strength of the magnetic field, but this remains an open problem.

Consider a convex billiard table having the symmetry of an ellipse. Such a billiard has a 2-periodic orbit along the $y$-axis.

\begin{problem}
Is it possible to choose the domain so that the dynamics of the corresponding billiard map are locally (near the 2-periodic orbit) conjugated to the dynamics of the rigid rotation through some angle   $\alpha$?
\end{problem}

 V. Schastnyy and D. Treschev \cite{ST}  provide numerical evidence that this is possible if $\alpha/\pi$ is  a Diophantine number.

\subsubsection{Geometry of caustics, invariant surfaces, and commuting billiard maps}

 A $p/q$-rational caustic corresponds to an invariant circle of the billiard map that consists of periodic points with the rotation number $p/q$. For example, a curve of constant width has a $1/2$ caustic, the envelop of the 2-periodic, back and forth, billiard orbits. A Birkhoff integrable billiard has a family of $p/q$-rational caustics for all sufficiently small values of $p/q$ -- this follows from the Arnold-Liouville theorem. An ellipse possesses $p/q$-rational caustics for all  $p/q< 1/2$.

 Kaloshin and K. Zhang \cite{KaZh} proved that, in an appropriate topology, the set of Birkhoff billiards that possess a rational caustic is dense in the space of Birkhoff billiards.

Baryshnikov and Zharnitsky  \cite{BZ} and Landsberg \cite{Lan} used methods of exterior differential systems to construct infinite-dimensional families of billiards having a $p/q$-rational caustic for a fixed value of $p/q$; see \cite{GenT} for outer billiards.

\begin{question}
Are there plane billiards, other than ellipses, that possess rational caustics with two different values of the rotation numbers?  Same question for outer billiards.
\end{question}

Kaloshin and J. Zhang \cite{KalZ} studied the case of $1/2$ and $1/3$ caustics. They proved that, for a certain class of analytic deformations of a circle, the coexistence of such caustics implies that the deformation is trivial.

Here is a weaker version of the previous problem. Let $q_0 \geq 2$ be an integer. A billiard table is called $q_0$-rationally integrable is it possesses $p/q$-rational caustics for all $0< p/q < 1/q_0$ (see \cite{Huang}).

 \begin{con}  \label{ratint} {\rm (V. Kaloshin)} {\it
A $q_0$-rationally integrable billiard is elliptic.}
 \end{con}

The next conjecture is motivated by a theorem of N. Innami's: if there is a sequence of convex caustics with rotation numbers tending to $1/2$, then the billiard table is an ellipse \cite{Inn} (see \cite{AB} for a simpler proof). An invariant circle of the billiard map is called {\it rotational} if it consists of periodic points.

 \begin{con}  \label{Inlim} {\rm (M. Bialy)} {\it
Assume that a plane billiard has a family of rotational invariant circles with all the rotation numbers in the interval $(0,1/2)$. Then the billiard table is an ellipse.
 }
 \end{con}

 A variation on Innami's result is the following question, suggested by V. Kaloshin:

 \begin{question}
Suppose that a billiard table has a sequence of convex caustics with rotation numbers converging to some number $\omega \in (0,1/2)$. Does it imply that the table is an ellipse?
\end{question}
The previously discussed rotation numbers, $0$ and $1/2$, are clearly special.

One may generalize the notion of caustic for higher-dimensional billiards as a hypersurface such that if a segment of the billiard trajectory is tangent to it, then the reflected segment is again tangent to this hypersurface. However, this notion is very restrictive: if the dimension of the Euclidean space is at least three, the only convex billiard table that possess caustics are ellipsoids \cite{Be,Gru}.

 \begin{question} {\rm (M. Bialy)}
Are there multi-dimensional convex billiards, other than ellipsoids, having invariant hypersurfaces in the phase space?
 Same question for multi-dimensional outer billiards.
\end{question}

 A common feature of families of integrable maps is that they commute (Bianchi permutability).

 Consider two confocal ellipses: the billiard reflections define two transformations on the set of oriented lines that intersect both ellipses. These transformations commute, see \cite{Tabachnikov2005}. The converse was recently proved by A. Glutsyuk \cite{Glutsyuk171}, as a consequence of his classification of 4-reflective complex planar billiards. For outer billiards, an analogous theorem is that if the outer billiard transformations on two nested curves commute then the curves are concentric homothetic ellipses \cite{Tabachnikov94}. It would be interesting to find a proof of Glutsyuk's result via a more geometric approach of \cite{Tabachnikov94}.

In higher dimensions (greater than two), the billiard transformations, induced by two confocal ellipsoids, commute as well (see, e.g., \cite{DR}).

\begin{problem}{\rm (A. Glutsyuk, S. Tabachnikov \cite{Tabachnikov2015})}
Given two nested closed convex hypersurfaces, assume that the billiard transformations on the set of oriented lines that intersect both hypersurfaces commute. Prove that the hypersurfaces are confocal ellipsoids.
\end{problem}

 The above problem also has analogs in the elliptic and hyperbolic geometries. Namely, one can define quadrics in the elliptic and hyperbolic geometries: a quadratic hypersurface of the unit sphere is, by definition, its intersection with a quadratic cone, and a similar definition applies to the hyperboloid model of the hyperbolic geometry. It is known that the billiard maps in quadrics are completely integrable in the geometries of constant curvature, see, e.g., \cite{Tabachnikov95}.

\subsection{Generalized pentagram maps}

We already mentioned the pentagram map. This map sends a (generic) polygon in the projective plane to a new polygon whose vertices are the intersection points of the consecutive short diagonals of the original polygon. The map extends to {\it twisted} $n$-gons, bi-infinite collections of point $(v_i)$ satisfying $v_{i+n} = M(v_i)$ for all $i$, where $M$ is a projective transformation. The pentagram map also commutes with projective transformations (since it only involves points and lines, and the operations of connecting pairs of points by lines and intersecting lines), hence the pentagram map is defined on the moduli space of projective equivalence classes of polygons.

Numerous generalizations of the pentagram map were introduced and studied, in particular, in \cite{KS1,KS2,KS3,MB}. Following \cite{KS2}, we describe one class of such generalized pentagram maps and present some problems and conjectures.

Let $(v_i)$ be a twisted polygon in $\mathbb {RP}^d$. Fix two $(d-1)$-tuples of integers
$$
I = (i_1,\ldots,i_{d-1}),\ J= (j_1,\ldots,j_{d-1}).
$$
Define the hyperplanes
$$
P_k = {\rm Span} (v_k, v_{k+i_1}, v_{k+i_1+i_2},\ldots, v_{k+i_1+\ldots+i_{d-1}}),
$$
and consider the new point
$$
T_{I,J} (v_k) = P_k \cap P_{k+j_1} \cap P_{k+j_1+j_2} \cap \ldots \cap P_{k+j_1+\ldots+j_{d-1}}.
$$
This formula defines the generalized pentagram map $T_{I,J}$. In this notation, the usual pentagram map is $T_{(2),(1)}$.

For some of these maps, for example, $T_{(2,2),(1,1)}, T_{(2,1),(1,1)}$, and $T_{(3,1),(1,1)}$, complete integrability has been established \cite{KS1,KS3,MB}. Other cases were studied numerically in \cite{KS2}.

This study was based on the height criterion \cite{Hal,Gram}. The height of a rational number $p/q$, given in the lowest terms, is $\max (|p|,|q|)$. In appropriate coordinates, $T_{I,J}$ is a rational map defined over $\mathbb Q$.
One considers the growth of the maximal heights of the coordinates  under iterations of the map. For integrable maps, one expects a slower growth, and for non-integrable ones -- a faster one. More specifically, a linear growth of the double logarithm of this height  is a numerical indication of non-integrability.

In particular, these experiments suggested that the maps
$$
T_{(2,2),(1,2)}, T_{(1,2),(1,2)}, T_{(1,3),(1,3)}, T_{(2,3),(2,3)},
$$
are integrable, whereas the maps
$$
T_{(1,2),(3,1)}, T_{(1,2),(1,3)}, T_{(2,3),(1,1)}, T_{(2,4),(1,1)}, T_{(3,3),(1,1)}
$$
are not integrable.

 \begin{con}  \label{hpent} {\rm (B. Khesin, F. Soloviev)}
The maps $T_{I,I}$ are integrable for all multi-indices $I$.

 \end{con}

 Khesin and Soloviev also conjectured that all the maps $T_{I,J}$ (and even a more general type of maps, called universal pentagram maps) are Hamiltonian, that is, preserve a certain non-trivial Poisson structure. See the recent paper by A. Izosimov on this subject \cite{IzosimovPent}. Some of the maps considered by Khesin and Soloviev are included in the family of  transformations described in \cite{GP}; they are examples of integrable cluster dynamics.

\subsection{Noncommutative integrable maps}

This section concerns dynamical systems, with continuous or discrete time, on associative but not commutative algebras, so that the variables involved do not commute; the reader with an aversion to free non-commuting variables may think of matrices of an arbitrary size: see, e.g., \cite{OlSo,MiSo}. In this case, one speaks about {\it noncommutative integrability}, although this term has another, quite different, meaning as well (when the integrals {form} a non-Abelian Lie algebra with respect to the Poisson bracket).

The overarching problem is to develop a comprehensive theory that would include noncommutative counterpart to the definition of integrability as given in Section \ref{defint}. As of now, the geometric aspects of the noncommutative theory are poorly understood, and one mostly deals with the formal, algebraic part of the story. In what follows, we  work with noncommutative Laurent polynomials or power series in several variables.

Here is an example of a 2-dimensional noncommutative integrable map that illustrates the specifics and peculiarities of the noncommutative situation. This example was discovered by S. Duzhin and M. Kontsevich \cite{Konts} and studied in \cite{Art,WE}:
$$
T: (x,y) \mapsto (xyx^{-1}, (1+y^{-1})x^{-1}).
$$
Here $x$ and $y$ are free non-commuting variables and the map can be understood as a substitution in the ring of Laurent polynomials $\bbC\langle x,y,x^{-1},y^{-1}\rangle$.

The map $T$ preserves the commutator $xyx^{-1}y^{-1}$. In the commutative case, when $xy=yx$, the map becomes
$$
T: (x,y) \mapsto (y, (1+y^{-1})x^{-1}),
$$
but the integral $xyx^{-1}y^{-1}$ degenerates to a constant.
On the other hand, in the commutative case, one has an integral
$$
I=x+y+x^{-1}+y^{-1}+x^{-1}y^{-1},
$$
whose level curves are elliptic curves, and the map has an invariant area form $(xy)^{-1} dx \wedge dy$, so it is completely integrable in the sense of Section \ref{defint}.

In the noncommutative case, the Laurent polynomial $I$ is transformed by $T$ as follows $I \mapsto xIx^{-1}$. If the variables represent matrices (of an arbitrary size), then $\tr I^k$ is an integral of the transformation $T$ for every positive integer $k$. In the general case, when $x,y$ are formal variables, one defines trace on the algebra of Laurent polynomials $\bbC\langle x,y,x^{-1},y^{-1}\rangle$ as an equivalence class of the cyclic permutations of factors of a polynomial.

On the other hand, in the intermediate (``quantum") case, when $xy=qyx$, with $q$ commuting with everything, the map $T$ has the integral
$$
I_q=x+qy+qx^{-1}+y^{-1}+x^{-1}y^{-1},
$$
whose $q \to 1$ limit  is $I$.

The map $T$ is a discrete symmetry of the system of differential equations
$$
\dot x = xy-xy^{-1}-y^{-1},\ \dot y = - yx+yx^{-1}+x^{-1}.
$$
Once again, we consider these differential equations in the formal sense, as Laurent series whose coefficients depend on a parameter, and we are not concerned with the questions of existence and uniqueness of solutions.

This system of ODEs admits a Lax presentation with a spectral parameter (that commutes with everything); the entries of the Lax matrix are noncommutative Laurent polynomials in $x$ and $y$, and the conjugacy class of the Lax matrix does not change under the map $T$.

The map $T$ belongs to a family of maps
$$
T_k: (x,y) \mapsto (xyx^{-1}, (1+y^{k})x^{-1}).
$$
For $k \ge 1$, these maps satisfy the {\it Laurent phenomenon}: both components of all the iterations belong to the noncommutative ring of Laurent polynomials in $x,y$ (for $k\ge 3$, the maps $T_k$ are not integrable).

This theorem, also conjectured by S. Duzhin and M. Kontsevich \cite{Konts}, was first proved in \cite{Usn} and, independently, in \cite{DiKe1}, and later, by different methods, in \cite{BeRe,Lee} (these papers also tackle different variations of the S. Duzhin and M. Kontsevich conjectures). In the commutative setting, the Laurent phenomenon is well known in the theory of cluster algebras and is related with integrability of cluster dynamical systems.

We now switch gears and discuss an open problem concerning noncommutative versions of the pentagram
maps described in the previous section. These maps belong to the realm of projective geometry, so one needs a noncommutative version of projective geometry and, in particular, of the cross-ratio of a quadruple of points on the projective line.

A fundamental notion of noncommutative algebra is quasi-determinant \cite{GeRe,GeRe1,GGRW}, a noncommutative analog of the determinant (quasi-determinants play a role in the study of noncommutative integrability, see \cite{EGR,DiKe2}). Using quasi-determinants, V. Retakh \cite{Ret} defined the cross-ratio and described its properties. Let us present the relevant definitions.

Let $A$ be a $2\times n$ matrix whose entries $a_{1i}$ and $a_{2i}$, $i=1,\ldots,n$, are elements of an associative ring $R$ over a field.
For $1 \leq i,j,k \leq n$, define the quasi-Pl\"ucker coordinates of the matrix $A$
$$
q_{ij}^k (A)= (a_{1i} - a_{1k} a_{2k}^{-1} a_{2i})^{-1} (a_{1j} - a_{1k} a_{2k}^{-1} a_{2j}),
$$
where we assume that all the inversions are defined (quasi-Pl\"ucker coordinates can be defined for $k \times n$ matrices as well).

Let $x,y,z,t$ be four 2-vectors with components in $R$. These vectors form the columns of a $2\times 4$ matrix $A$, and one defines the cross-ratio
$$
[x,y,z,t] = q_{zt}^y (A) q_{tz}^x (A).
$$
In the commutative case, assuming that the second components of the  vectors equal 1, one obtains a familiar definition of the cross-ratio (or, better said, one of its six definitions):
$$
[x,y,z,t] = \frac{(z-x)(t-y)}{(z-y)(t-x)}.
$$
The non-commutative cross-ratio shares some properties with its commutative counterpart. For example,
$$
[G x \lambda_1, G y \lambda_2, G z \lambda_3, G t \lambda_4] = \lambda_3^{-1} [x,y,z,t] \lambda_3,
$$
where $G$ is a $2\times 2$ matrix over $R$, and $\lambda_i \in R,\ i=1,2,3,4$.

Having a definition of noncommutative cross-ratio, one may define and study noncommutative analogs of the pentagram maps.
To be concrete, consider the following 1-dimensional version of the pentagram map, introduced in  \cite{GSTV} and called the {\it leapfrog map}.

Let $S^-=(S_1^-,S_2^-,\ldots)$ and $S^-=(S_1,S_2,\ldots)$ be a pair of $n$-gons in the projective line. The leapfrog map $(S^-,S)\mapsto (S,S^+)$ is defined as follows. For every $i$, there exists a unique projective involution that sends the triple $(S_{i-1},S_i,S_{i+1})$ to the triple  $(S_{i+1},S_i,S_{i-1})$; this transformation sends $S_i^-$ to $S_i^+$ (that is, the point $S_i^-$ ``jumps" over $S_i$, in the projective metric on the segment $[S_{i-1},S_{i+1}]$, and becomes $S_i^+$ -- thus the name of the map). In other words, one has
$[S_{i-1},S_i,S_{i+1},S_i^-]=[S_{i+1},S_i,S_{i-1},S_i^+]$.

The theory of noncommutative cross-ratios implies that there exists a projective transformation
$
\left(S_{i-1}, S_i, S_{i+1}, S_i^{-}\right) \mapsto \left(S_{i+1} , S_i , S_{i-1} , S_i^{+}\right)
$
 if and only if
\begin{equation*}
\begin{split}
&\left(S_{i+1}-S_i \right)^{-1}\left(S_i^{-} - S_i\right)\left(S_i^{-} - S_{i-1}\right)^{-1}\left(S_{i+1} - S_{i-1}\right)=\\
r^{-1}&\left(S_{i-1} - S_i\right)^{-1}\left(S_i^{+}-S_i\right)\left(S_i^{+} - S_{i+1}\right)^{-1}\left(S_{i-1} - S_{i+1}\right)r
\end{split}
\end{equation*}
where $r \in R$ is  invertible. In the commutative case, this becomes the formula defining the leapfrog map.

It was proved in \cite{GSTV} that the leapfrog map is completely integrable: it has an invariant Poisson structure and a complete collection of Poisson commuting integrals; in the appendix to \cite{GSTV}, A. Izosimov constructed a tri-Hamiltonian structure for this map.

\begin{problem}{\rm (V. Retakh)
Establish complete integrability of the noncommutative version of the leapfrog map. Define noncommutative versions of the pentagram map and its higher-dimensional analogs and establish their complete integrability.
}
\end{problem}

Let us add that a version of the pentagram map was defined and studied for polygons in Grassmannians \cite{FeMa}. It is a further challenge to define a noncommutative version of these maps and to investigate their integrability.


\section{Bi-Hamiltonian integrable systems}\label{bi}


Let $M$ be a smooth manifold. A \textit{Poisson bracket} on $M$ is a bilinear (over $\bbR$) skew-symmetric operation $\{~{,}~\}$
on the space of smooth functions $C^\infty(M)$   which satisfies the Leibniz rule
$$
\{fg,h\}=f\{g,h\}+\{f,h\}g
$$
and the Jacobi identity
$$
\{f,\{g,h\}\}+\{h,\{f,g\}\}+\{g,\{h,f\}\}=0
$$
for any $f,g,h\in C^\infty(M)$.
This operation turns $C^\infty(M)$ into an infinite-dimensional Lie algebra.

In local coordinates,  a Poisson bracket can be written as
$$
\{f, g\} (x) = \sum_{i.j}^n A^{ij}(x) \frac{\partial f}{\partial x_i}\frac{\partial g}{\partial x_j},
$$
where $A=\bigr( A^{ij}(x)\bigl)$ is the corresponding Poisson tensor. If the rank of $A$  is smaller than the dimension of $M$, the  basic  definition of a finite-dimensional integrable systems given in Introduction needs a small modification. Namely, we will assume that the number of independent commuting first integrals (including the Hamiltonian $H$) equals
$$\frac{1}{2}(\dim M + \operatorname{corank} A )$$
where the corank of $A$ is taken at a generic point.  This condition is equivalent to the fact that the subspace spanned by the differentials of the first integrals in the cotangent space $T^*_xM$ at a generic point $x\in M$ is maximal isotropic with respect to the Poisson tensor $A$  (which is exactly the same condition as what we use in the symplectic case).

One of the most effective methods for constructing and studying integrable systems is based on the notion of compatible Poisson structures, see for example \cite{Magri78, Magri-Morosi}.
Recall that two Poisson brackets $\{~{,}~\}_1$ and $\{~{,}~\}_2$ on  $M$ are called \textit{compatible}
if any linear combination of them is again a Poisson bracket. The set  of all non-trivial linear combinations of two compatible Poisson brackets $\{~{,}~\}_1$ and $\{~{,}~\}_2$ is called a \textit{Poisson pencil}.

\subsection{Bi-Poisson vector spaces}

At a fixed point $x\in M$,  this pencil can be understood as a pencil of skew-symmetric bilinear forms $J=\{A_\lambda = A + \lambda B\}_{\lambda \in \bar{\mathbb C}}$  on the tangent space $T^*_xM$  (we may consider these forms  up to proportionality and formally set $A_\infty = B$). To each pencil $J$, we can naturally assign its algebraic type.  To explain what is meant by the algebraic type of $J$, we recall the Jordan-Kronecker decomposition theorem  (see, e.g., \cite{BolsZhang, tom})  that provides a simultaneous canonical form for a pair of skew-symmetric forms.

Consider a complex vector space $V$ with a pair of bilinear forms $A, B: V\times V \to \mathbb C$.  For simplicity we will assume that $B$ is regular in the pencil $J=\{A_\lambda = A + \lambda B\}$, i.e., the rank of $B$ is maximal within this pencil.

\begin{Th} \label{lin2}
By choosing an appropriate basis,  $A$ and $B$ can be simultaneously  reduced to a block-diagonal form:
\begin{equation}
\label{JKdecom}
A   \mapsto \begin{pmatrix}   A_{1} & &    \\  & \ddots &  \\  & & A_s    \end{pmatrix}, \quad
B   \mapsto \begin{pmatrix}   B_{1} & &   \\  & \ddots &  \\  & & B_s    \end{pmatrix}
\end{equation}
where the corresponding pairs of blocks have one of the following two possible types{\/\rm:}
\begin{align*}
&1) \mbox{ Jordan type $\mu_i$-block}: \quad
A_i = \begin{pmatrix}\hphantom-0& \mathrm{Id}\\- \mathrm{Id}&0\end{pmatrix},\quad
B_i = \begin{pmatrix}\hphantom-0&J(\mu_i)\\-J^\top(\mu_i)&0\end{pmatrix},\\
\intertext{where $J(\mu_i)$ is a Jordan block of size $k_i\times k_i$ with eigenvalue $\mu_i$ and $\mathrm{Id}$ is the identity matrix of the same size{\/\rm;}}
&2)\mbox{ Kronecker block}:\quad
A_i=\begin{pmatrix}0&D\\-D^\top&0\end{pmatrix},\quad
B_i=\begin{pmatrix}0&D'\\-{D'}^\top&0\end{pmatrix},
\end{align*}
where $D$ and $D'$ are $k_i\times(k_i+1)$--matrices of the form:
$$
D = \begin{pmatrix}
  1 & 0  & & \\
      &\ddots & \ddots & \\
     &           &  1 & 0
     \end{pmatrix}, \quad D' =  \begin{pmatrix}
  0 & 1  & & \\
      &\ddots & \ddots & \\
     &           &  0 & 1
     \end{pmatrix}
$$
\end{Th}

By the {\it algebraic type} of the pencil $J=\{A_\lambda = A + \lambda B\}$ we mean the type of the above decomposition that includes the number and sizes of  Kronecker blocks and Jordan blocks, separately for each characteristic number $\mu_i$ (ignoring, however, the specific values of $\mu_i$'s). It follows from this definition that, in each dimension, there are only finitely many different algebraic types.

Clearly, the algebraic type of a Poisson pencil may depend on a point $x\in M$  but, in the real analytic case,  there always exists an open everywhere dense subset $U\subset M$ of generic points with the same algebraic type.   We will refer to it as the algebraic type of the Poisson pencil on $M$.

It is well known that many properties of a bi-Hamiltonian system essentially depend on the algebraic type of the underlying Poisson pencil. For instance, if this pencil is Kronecker (i.e., there are only Kronecker blocks in the Jordan-Kronecker decomposition), then the integrals of any related bi-Hamiltonian system are Casimir functions of the brackets forming the pencils, whereas in the case of Jordan (symplectic) pencils the integrals are the eigenvalues of the so-called recursion operator (equivalently, characteristic numbers of the Jordan-Kronecker decomposition). For Poisson pencils of a ``mixed'' algebraic type, the situation becomes more complicated. In particular, there might be no canonical choice of integrals for associated bi-Hamiltonian systems.  One of the most natural problems in this context would be to develop a theory of bi-Poisson vector spaces, i.e., spaces endowed with a pencil of skew symmetric forms. Such a theory would serve as a natural analog of symplectic linear algebra in view of the role the latter plays in symplectic geometry and the theory of Hamiltonian systems.  Some of results in this direction can be found in \cite{BolsZhang, Pumei}, but still  there is a number of further questions to be clarified. Below we list some of them.

The automorphism group of a pencil $J=\{A_\lambda = A + \lambda B\}$ is an algebraic linear group defined by
\begin{equation} \label{automorph}
\mathrm{Aut} (V,J)=\{\phi\in\mathrm{End}(V)\mid A_\lambda(\phi(\xi),\phi(\eta))=A_\lambda(\xi,\eta) \ \mbox{for all } A_\lambda \in J\}.
\end{equation}

Obviously,  the structure of this group essentially depends on the algebraic type of $J$ (see \cite{Pumei}).

\begin{problem}
Consider the action of $ \mathrm{Aut} (V,J)$ on  $V$.  Describe the partition of $V$ into $ \mathrm{Aut} (V,J)$-orbits.  More generally,  describe the action of $ \mathrm{Aut} (V,J)$ on the set of all $k$-dimensional subspaces $U\subset V$ (orbits, invariants, fixed points).  Notice that fixed points of this action (i.e., $ \mathrm{Aut} (V,J)$-invariant subspaces)  are important, as they correspond to well-defined (co)distributions in the context of Poisson pencils.
\end{problem}

A subspace $L\subset V$ is said to be {\it bi-Lagrangian} if it is simultaneously isotropic with respect to both  $A$ and $B$ and has maximal possible dimension  (which is $\frac{1}{2}(\dim V + \corank J)$). The set $\mathrm{LG}(V,J)$ of all bi-Lagrangian subspaces is called the {\it bi-Lagrangian Grassmannian}.  A detailed description of $\mathrm{LG}(V,J)$ is an important open problem, see the discussion in \cite{RosScho}. It is important that bi-Lagrangian subspaces always exist, so that  $\mathrm{LG}(V,J)$ is never empty.  In general, this is a projective algebraic variety with rather non-trivial properties, which essentially depend on the algebraic type of the pencil $J$.  In particular, there might exists bi-Lagrangian subspaces of different algebraic types. In other words, $\mathrm{LG}(V,J)$ may consist of several orbits of the automorphism group $ \mathrm{Aut} (V,J)$.

\begin{problem}
Find necessary and sufficient conditions for the bi-Lagrangian Grassmannian $\mathrm{LG}(V,J)$ to be a smooth algebraic variety.
Describe the partition of   $\mathrm{LG}(V,J)$  into $\mathrm{Aut} (V,J)$-orbits.
\end{problem}

A ``non-linear analog''  of the existence of bi-Lagrangian subspaces would be the  existence of  complete subalgebras  $\mathcal F\subset C^\infty(M)$ consisting of functions in involution with respect to all the brackets from a given Poisson pencil. By completeness in this case we mean that the subspace generated by the differentials of functions $f\in\mathcal F$ in the cotangent space $T^*_x M$ is maximal isotropic. It is natural to refer to such a subalgebra (or to a collection of its generators) as a {\it bi-integrable} system.

In the case of Kronecker pencils, such a subalgebra $\mathcal F$ is generated by the (local) Casimir functions of  the brackets in the pencil.  It is well known that $\mathcal F$ so obtained is complete and commutative with respect to all the brackets in the pencil  (\cite{Bols91, gelzak1}). Similar results in the symplectic case are due to P.\,Olver~\cite{olver}  and H.\,Turiel~\cite{turiel}.  In the mixed case  (when the Jordan-Kronecker decomposition of a pencil includes both Kronecker and Jordan blocks),  bi-integrable systems do exist locally in a  neighbourhood of a generic point. However,  to the best of our knowledge, this result has never been published.

In view of the above discussion on the partition of $\mathrm{LG}(V,J)$ into different orbits,  it would be natural to distinguish bi-integrable systems of different algebraic types. More precisely,  let  $f_1, \dots, f_k$ be a complete set of functions in bi-involution.   We restrict our considerations to a small neighbourhood $U$ of a generic point $x\in M$ such that for each $y\in U$ the Jordan-Kronecker type of the pencil defined on $T^*_y M$ by a pair of compatible Poisson brackets is the same. If we think of $T^*_y M$ as a bi-Poisson vector space, then the subspace $L_y=\mathrm{span}\{df_1(y), \dots, df_k(y)\} \subset  T^*_y M$ is bi-Lagrangian and therefore can be characterised by its algebraic type (i.e., by the type of the orbits of the automorphism group action in the bi-Lagrangian Grassmannian to which $L$ belongs).  If this type remains the same for all points $y\in U$, we refer to it as the {\it algebraic type} of
the bi-integrable system $f_1, \dots, f_k$.

\begin{problem}
Do bi-integrable systems exist for each algebraic type?
\end{problem}

This question is local in the sense that bi-integrable systems in question ``live'' in a small  neighbourhood  of a generic point $x\in M$  where the algebraic type of a given Poisson pencil is locally constant. However,  for applications, we need smooth commuting  functions that are globally defined on the whole manifold $M$.  For two constant brackets $\{~,~\}_1$ and $\{~,~\}_2$ on a vector space $V$, bi-integrable systems can be easily constructed by means of bi-Poisson linear algebra: each bi-Lagrangian subspace can be interpreted as a bi-integrable Hamiltonian system with linear integrals.  The first non-trivial example is the so-called ``argument shift'' pencil (see A.S.\,Mischenko and A.T.\,Fomenko~\cite{MF78-Izv}), generated by a linear and a constant brackets.  We recall this construction.

\subsection{Argument shift method and Jordan-Kronecker invariants of Lie algebras}\label{argshift}

Let $\goth g^*$ be the dual space of a finite-dimensional Lie algebra $\goth g$.  It is well known that~$\goth g^*$ possesses two natural compatible Poisson brackets.  The first one is the standard linear Lie-Poisson bracket
\begin{equation}
\label{LiePoisson}
\{ f, g \} (x) =\langle x ,[d f(x), d g(x)] \rangle,
\end{equation}
and the second one is a constant bracket given by
\begin{equation}
\label{abracket}
\{ f, g \}_a (x) = \langle a ,[d f(x), d g(x)] \rangle,
\end{equation}
where $a\in\goth g^*$ is a fixed element.
Here we assume  $a$ to be regular, although formula \eqref{abracket} makes sense for an arbitrary $a$.

There are many examples of finite-dimensional Lie algebras for which we can construct a complete family of polynomials in bi-involution with respect to the  brackets \eqref{LiePoisson} and \eqref{abracket}.  In particular, the famous theorem by A.S.\,Mischenko and A.T.\,Fomenko~\cite{MF78-Izv} states that such a family exists for every semisimple Lie algebra $\goth g$ and consists of shifts of $\Ad^*$-invariant polynomials, i.e., function of the form $f_i(x + \lambda a)$, where  $\lambda\in \R$ and $f_1,\dots, f_r$ are generators of the algebra of (co)-adjoint invariants.

However, it is still an open question whether  one can
construct  such a family  for every finite-dimensional Lie algebra.  In many examples (see \cite{BolsZhang}),  the answer turns out to be positive, leading us to  the following Generalised Argument Shift conjecture which was first stated in \cite{BolsZhang}  and then mentioned and discussed in a number of papers \cite{RosScho, BolsIzosTson,  IzosimovGAS, bolsMFconj}:

\begin{con}  Let $\goth g$ be  a finite-dimensional Lie algebra. Then for every regular element $a\in \goth g^*$,   there exists a complete family $\mathcal G_a$ of  polynomials in bi-involution, i.e., in  involution with respect to the two brackets
$\{\,,\,\}$ and $\{\,,\,\}_a$. \end{con}

It is expected that  the solution to  this conjecture  will essentially rely on the algebraic type of the pencil  $\{~,~\}_\lambda = \{~,~\} + \lambda\{~,~\}_a$ or, equivalently, of the pencil of the skew-symmetric forms $\cA_{x+\lambda a}$, $\lambda\in \mathbb C$,  on $\goth g$ for a generic pair $(x,a)\in \goth g^*\times \goth g^*$, where
$$
\cA_{x+\lambda a}(\xi,\eta)=\langle x+\lambda a, [\xi,\eta] \rangle, \quad \xi,\eta\in\Fg.
$$
This algebraic type can be considered as an invariant of the Lie algebra $\Fg$, called {\it  Jordan--Kronecker invariant} (see \cite{BolsZhang}).

\begin{problem}
Compute the Jordan--Kronecker invariants for  the most interesting classes of Lie algebras, and particularly  for the following
\newline
{\/\rm(}a{\/\rm)} semidirect sums $\Fg+_{\rho}V$, where $\rho: \Fg \rightarrow \mathrm{End}(V)$ is a representation of a semisimple Lie algebra $\Fg$ and $V$ is assumed to be commutative;
\newline
{\/\rm(}b{\/\rm)} Borel subalgebras of simple Lie algebras;
\newline
{\/\rm(}c{\/\rm)} parabolic subalgebras of simple Lie algebras;
\newline
{\/\rm(}d{\/\rm)} the centralisers of singular elements $a\in\Fg$, where $\Fg$ is simple;
\newline
{\/\rm(}e{\/\rm)} Lie algebras of small dimensions.
\end{problem}

For Lie algebras of dimension $\le 5$, the Jordan--Kronecker invariants were computed by P.\,Zhang. Her results, as well as some other examples of computing Jordan--Kronecker invariants, can be found in the arXiv version of~\cite{BolsZhang}.  For some semidirect sums,  Jordan-Kronecker invariants have been recently computed by K.\,Vorushilov \cite{Vorushilov}.

\begin{question}\label{prob25}
Are there any restrictions on the algebraic type of the pencils $\cA_{x+\lambda a}$? Which algebraic types can be realised by means of an appropriately chosen Lie algebra{\/\rm?}
\end{question}

A recent observation by I. Kozlov (see details in \cite{BolsZhang})  shows that some restrictions do exist.  Roughly speaking, Kozlov has obtained the answer  for Lie algebras of pure Jordan and pure  Kronecker type.  The mixed case remains open because of possible non-trivial interaction between Kronecker and Jordan blocks.

To further develop these ideas, it would be interesting to study Poisson pencils generated by one quadratic and one linear brackets.  Compatible brackets of this kind appear, in particular, in the theory of Poisson-Lie groups.

\begin{problem} Study examples of  ``quadratic $+$ linear'' Poisson pencils.
Compute their algebraic types and construct complete families of polynomials in bi-involution for such pencils.
\end{problem}

It is a very general question whether or not a given symplectic  (Poisson) manifold $M$ admits an integrable system with good properties. In particular, if $M$ is a vector space endowed with a polynomial Poisson bracket, we are interested in integrable systems defined on $M$ by a collection of Poisson commuting polynomials.
For a linear Poisson bracket  \eqref{LiePoisson}, such a polynomial integrable system always exists. This statement, known as Mischenko-Fomenko conjecture, was proved by S. Sadetov \cite{sadetov} in 2004.   What can we say about polynomial Poisson brackets of higher degrees?

\begin{question}
Is it true that for any quadratic Poisson bracket  (defined on a vector space),  there exists a polynomial integrable system?
\end{question}

\subsection{Nijenhuis operators: singular points and global properties}

Let $L= \bigr(L^i_j(x)\bigl)$ be a $(1,1)$-tensor field (field of endomorphisms) on a smooth manifold $M$.

 We   say that $L$ is {\it a Nijenhuis operator} if its Nijenhuis torsion $N_L$ identically vanishes on $M$, i.e.,  $$
N_L(v, w) = L^2 [v, w] + [Lv, Lw] - L [Lv, w] - L [v, Lw] = 0,
$$
for arbitrary vector fields $v, w$ on $M$.

It is known that Nijenhuis operators are closely related to compatible Poisson brackets and therefore naturally appear in the theory of finite dimensional integrable systems, see e.g. Magri {\it et al} \cite{mks,Magri-Morosi}. Indeed, if a Poisson structure  given by the tensor
$P^{ij}$  is compatible with the Poisson structure given by the symplectic form $\omega$, then $L^{i}_j:= P^{is}\omega_{sj}$ is a Nijenhuis  operator. Nijenhuis tensors appear also in other topics in the theory of integrable systems, see e.g. \cite{benenti, Ibort}, in differential geometry, see e.g. \cite{bartelmeo, splitting, local},   so investigation of Nijenhuis tensors is expected to be generally important for the theory of integrable systems and other branches of mathematics;  in this section we collect open problems about it.

Let $L=L(x)$ be a $(1,1)$-tensor field. Then at each point $x\in M$, we can characterize $L(x)$ by its algebraic type (Segre characteristic), that is, the type of its Jordan normal form  (eigenvalues with multiplicities + sizes of Jordan blocks for each eigenvalue). The algebraic type of $L$ depends on $x$ and may vary from point to point.    We will say that $x\in M$ if {\it generic}, if the algebraic type of $L$ is locally constant at $x$  (i.e., is the same for all points from a certain  neighbourhood  $U(x)$).  Otherwise,  $x$ is called {\it singular}. One can show that generic points form an open  dense subset of $M$.

Almost all results known in the literature concern local properties of Nijenhuis operators  at generic points, e.g.  for semisimple Nijenhuis operators a local description is due to Haantjes \cite{Haantjes}.
   However,  Nijenhuis operators also occur  in  problems which require a global treatment and analysis of singular points.
Indeed, in the theory of bi-Hamiltonian systems, the Nijenhuis operators appear as recursion operators $R$  (see  \cite{Magri-Morosi, mks}).  The first integrals (for an important class of integrable systems)  are the coefficients of the characteristic polynomial of  $R$. This means, in particular, that the structure of singularities of the corresponding Liouville foliations  is determined by the local structure of singular points of $R$.  Analysing them, we may find topological obstruction to bi-Hamiltonisation, which would be an interesting result.  Philosophically,  we want to say that the  singularities of bi-integrable systems are determined by the singularities of the underlying structure  (in our case, this is just the recursion operator).  We also notice that integrability of finite-dimensional Hamiltonian systems is essentially a global phenomenon, so that it might be important to understand global properties of recursion operators.

Similarly,  in differential geometry,  Nijenhuis operators play an important role in the theory of projectively equivalent (pseudo-)Riemannian metrics, where understanding of singular points of related Nijenhuis tensors
led to the proof  of  the projective Lichnerowicz conjecture in the Riemannian and Lorentzian signature \cite{lichnerowicz, BMR}.

In view of the above discussion, we find the following directions of research  very important:
\begin{itemize}
\item Studying global properties of Nijenhuis operators  and topological obstructions to the existence of Nienhuis operators (of a certain type) on compact manifolds.

\item Studying  local behaviour of Nijenhuis operators at singular points.
\end{itemize}

Speaking of topological obstructions to the existence of Nijenhuis opertors on compact manifolds,  we should take into account two simple but important observations:
\begin{itemize}
\item[(i)]  The operators of the form $f(x)\cdot \mathrm{Id}$ are Nijenhuis for any smooth function $f:M \to \R$.  Such operators exist on any manifold but this example is too trivial.
\item[(ii)]  In $\R^n$, there exist  $C^\infty$-smooth Nijenhuis operators with compact support and such that, in a certain domain, they take the standard diagonal form  $L(x) = \mathrm{diag} ( \lambda_1(x_1), \lambda_2(x_2), \dots , \lambda_n(x_n))$ with distinct non-constant real eigenvalues $\lambda_i(x_i)$.   Since the support of such an operator is contained in a certain ball,  we can construct a Nijenhuis operator with similar properties on any (compact) manifold.
\end{itemize}

This leads to the following  natural questions.

We will say that a linear operator $L: V\to V$ is {\it algebraically regular} if the dimension of the orbit of $L$ under the canonical $GL(V)$-action on the space of operators is maximal (equivalently, the dimension of the  centraliser of $L$ is minimal, i.e., equals  $n=\dim V$).  In particular,  diagonalisable operators with simple spectrum are regular.   For non-diagonalisable operators,  the algebraic regularity condition means that geometric multiplicity of each eigenvalue of $L$ is one.

\begin{problem}
Describe closed manifolds $M$ which admit Nijenhuis operators $L(x)$  that are algebraically regular at each point $x\in M$. 
\end{problem}

\begin{problem}\label{prob:compsupp} 
 Let us fix a certain algebraic type of a linear operator, i.e., its Segre characteristic (see above).   Does there exist a Nijenhuis operator $L$ in $\R^n$ with the following two properties:
\begin{enumerate}
\item  $L$ has compact support;
\item in a certain domain $U\subset \R^n$,  the algebraic type of $L$ does not change and coincides with the given one?
\end{enumerate}
More generally,  find necessary and sufficient conditions, in terms of the Segre characteristic, for the existence of a Nijenhuis operator in $\R^n$ with the above properties.
\end{problem}

 Notice that for semisimple operators $L$ with real eigenvalues, the answer to the above question is positive.  However, if some of the eigenvalues of $L$ are complex, then  $L$ cannot have compact support.  Also,  one can easily construct an example of a Nijenhuis operator $L$ with compact support which has the form of a Jordan block, e.g., $L = \begin{pmatrix}  y & x \\ 0 & y  \end{pmatrix}$, on a certain disc.  However, even for single Jordan blocks of size $\ge 3$ the answer is unknown.

\begin{problem}
Construct real analytic examples of Nijenhuis operators on closed two-dimensional surfaces whose eigenvalues are real and generically distinct.
\end{problem}

Let $\chi (t) = \det (t\cdot \Id - L)=t^n +\sigma_1(x) t^{n-1} + \sigma_2(x) t^{n-2} + \dots + \sigma_n(x)$ be the characteristic polynomial of a Nijenhuis operator $L$. Then the following relation holds  \cite{splitting}:
\begin{equation}
\label{eq10.1}
J L  = \begin{pmatrix}  \!\!\! -\sigma_1 & 1 & & & \\
 \!\!\! -\sigma_2 & 0 & 1 &  & \\
\vdots &  & \ddots & \ddots &  \\
-\sigma_{n-1} & & & 0 & 1 \\
 \!\!\! - \sigma_n & & & & 0   \end{pmatrix} J, \qquad\mbox{where } J =
 \begin{pmatrix}
 \frac{\partial \sigma_1}{\partial x_1} & \dots &  \frac{\partial \sigma_1}{\partial x_1} \\
 \vdots & \ddots & \vdots \\
 \frac{\partial \sigma_n}{\partial x_1} & \dots &  \frac{\partial \sigma_n}{\partial x_1}
 \end{pmatrix}
 \end{equation}

This formula shows that,  in the case when the coefficients of the characteristic polynomial of $L$ are functionally independent,   a Nijenhuis operator $L$ can be explicitly reconstructed from $\sigma_1, \dots, \sigma_n$.   Notice that, by continuity, this is still true if $\sigma_1, \dots, \sigma_n$ are independent almost everywhere on $L$.  In particular, the singularities of $L$  will be completely determined by the singularities of the smooth map $\Phi=(\sigma_1, \dots, \sigma_n) : M \to \R^n$.  However,  the singularities of $\Phi$ must be very special.

We can illustrate this with a very simple example.  Let $\dim M=2$, then we only have two coefficients $\sigma_1$ and $\sigma_2$.   Consider a point where $d\sigma_1$ and $d\sigma_2$ are linearly dependent and $d\Phi$ has rank 1.  Assume, in addition, that $d\sigma_1\ne 0$ so that there exists a local coordinate system $(x,y)$ such that $\sigma_1 = x$  (and $\sigma_2=f(x,y)$).  Then \eqref{eq10.1} implies that, in coordinates $(x,y)$,  our Nijenhuis operator $L$ takes the form
$$
L =\frac{1}{f_y}  \begin{pmatrix}
f_y & 0 \\ -f_x & 1
\end{pmatrix}
\begin{pmatrix} - x & 1 \\ -f & 0 \end{pmatrix}
\begin{pmatrix}
1 & 0 \\ f_x & f_y
\end{pmatrix} =  \begin{pmatrix}
f_x - x & f_y \\
\frac{f_x(x-f_x) - f}{f_y} & -f_x
\end{pmatrix}.
$$
Thus, $\sigma_2=f(x,y)$ is admissible if and only if $\frac{f_x(x-f_x) - f}{f_y} $  is a smooth function.  In other words, to describe singular points of $L$ in this rather simple particular case we need to solve the following

\begin{problem} Describe all the functions $f(x,y)$ of two variables defined in a neighbourhood of $(0,0)\in\R^2$ such that 
\begin{itemize}
\item $f_y(0,0)\not\equiv 0$;
\item $f_y(0,0)=0$; 
\item $\dfrac{f_x(x-f_x) - f}{f_y} $ is locally smooth (or, more precisely, extends up to a locally smooth function).
\end{itemize}
\end{problem}

A more general version of this question  (which would immediately clarify the local structure of singular points for a wide class of Nijenhuis operators)  is as follows.

\begin{problem}
Consider a smooth map $\Phi=(\sigma_1,\dots,\sigma_n):  U(0) \to \R^n$, where $U(0)$ is a  neighbourhood  of the origin $0\in \R^n$.  We assume that $0$ is a singular point of $\Phi$, that is, $\rank d\Phi(0) < n$, but almost all points $x\in U(0)$ are regular. Describe those maps $\Phi$ for which all the components of the Nijenhuis operator $L(x)$ defined by \eqref{eq10.1} are smooth  functions  (notice that $L$ is well defined and smooth at regular points of $\Phi$, and we are interested in necessary and sufficient conditions for $L$ to be smoothly extendable onto the whole  neighbourhood  $U(0)$).
\end{problem}

  The next portion of problems is related to   Nijenhuis operators all of whose components are linear functions in local coordinates $x_1, \dots, x_n$,   that is,
\begin{equation}
\label{eq:lN}
L ^i_j (x) = \sum_k b^i_{jk} x_k
\end{equation}

 The role of such operators in the theory of Nijenhuis operators is similar to the role of linear Poisson brackets in Poisson geometry.  Recall that linear Poisson brackets are in a natural one-to-one correspondence with finite-dimensional Lie algebras.  In the context of Nijenhuis operators, we have a similar statement (see, e.g., \cite{winterhalder}).
\begin{proposition}
An operator $L(x)$ with linear entries \eqref{eq:lN}  is Nijenhuis if and only if $b^i_{jk}$ are structure constants of a left-symmetric algebra.
\end{proposition}

Recall that an algebra $(A, *)$ is called {\it left-symmetric}  if the following identity holds:
$$
a * (b*c) - (a*b)* c = b*(a*c) - (b*a)*c,  \quad \mbox{for all } a,b,c\in A.
$$
The definition is due to Vinberg \cite{vinberg};  recent papers on left-symmetric algebras include \cite{burde,burde2}.

Only few results are known about  left-symmetric algebras. In particular, constructing new examples and classifying left-symmetric algebras of low dimensions would be a very interesting problem in the context of the theory of Nijenhuis operators.

For a given left-symmetric algebra with the structure constants $b^i_{jk}$, consider the coefficients $\sigma_1, \dots, \sigma_n$ of the characteristic polynomial of  $L(x)= \Bigl(L^i_j = \sum_k b^i_{jk} x_k\Bigr)$.  It is easy to see that $\sigma_k$ is a homogeneous polynomial of degree $k$.   Consider the class of left-symmetric algebras for which these polynomials are algebraically independent.  From \eqref{eq10.1} we conclude that the classification of such left-symmetric algebras is equivalent to that of collections of polynomials $\sigma_1, \dots, \sigma_n$, $\deg \sigma_k=k$, satisfying the following rather exceptional property.

Consider the Jacobian matrix $J = \left(   \frac{\partial (\sigma_1,\dots, \sigma_n)}{\partial(x_1,\dots, x_n)}\right)$.  It is easy to see that $D=\det J$ is a homogeneous polynomial of degree $\frac{n(n-1)}{2}$.   Consider the matrix
$$
L=J^{-1}  \begin{pmatrix}  -\sigma_1 & 1 & & \\
\vdots &   & \ddots & \\
-\sigma_{n-1} & & & 1\\
-\sigma_n & 0 & \dots & 0
\end{pmatrix} J, \quad\mbox{where } J = \left(   \frac{\partial (f_1,\dots, f_n)}{\partial(x_1,\dots, x_n)}\right)
$$
It is easy to check that the entries of this matrix are rational functions of the form $L_{ij}=\frac{Q_{ij}}{D}$ where $Q_{ij}$ are homogeneous polynomials of degree $\deg D + 1$. In other words, the entries of $L$ are homogeneous rational functions of degree one.  For some special polynomials $\sigma_1, \dots, \sigma_n$, a miracle happens:  all $Q_{ij}$ turn out to be divisible by $D$ and the entries of $L$ become linear functions in $x_1,\dots, x_n$.  Such $L$ are automatically Nijenhuis operators (equivalently, LSA's).

\begin{problem}
Describe/classify the collections of algebraically independent homogeneous polynomials $\sigma_1, \dots, \sigma_n$, $\deg \sigma_k = k$, in $n$ variables $x_1,\dots, x_n$  with the required property.
\end{problem}

Finally, one question revealing the relationship between singularities of bi-Hamiltonian systems and singular points of Nijenhuis operators.
In a recent paper \cite{BolsIzosCMP}, it was shown that the structure of singularities of bi-Hamiltonian systems related to Poisson pencils of Kronecker type is determined by the singularities of the pencil.  We expect that a similar principle works for Poisson pencils of sympletic type, i.e., when one of two compatible Poisson brackets is non-degenerate. In this case,  we have a well-defined recursion operator $R$ (which is automatically a Nijenhuis operator) and the commuting functions we are interested in are just the coefficients $\sigma_1,\dots,\sigma_n$ of the characteristic polynomial of $R$.  We assume that they are independent almost everywhere on $M$ and want to study singular points of the mapping
$\Phi=(\sigma_1,\dots,\sigma_n):  M \to \R^n$  playing the role of the momentum mapping for our bi-Hamiltonian system.

\begin{problem}   Describe the structure of singularities of $\Phi=(\sigma_1,\dots,\sigma_n):  M \to \R^n$ in terms of the singular points of the recursion operator $R$, their linearisations, and the corresponding left-symmetric algebras.   What are sufficient and/or necessary conditions for such singularities to be non-degenerate  (in the sense of Eliasson \cite{Eliasson})?  Notice that, in this case, the algebraic type of $R$ is rather special: at a generic point, $R$ is semisimple, but each of its eigenvalues has multiplicity 2.
\end{problem}


\section{Poisson geometry and action-angle variables}

Poisson manifolds constitute a natural generalisation of symplectic manifolds and their Hamiltonian dynamics generalises that of symplectic manifolds and adds singularities to its complexity  (see the beginning of Section \ref{bi} for basic definitions).

Let  $\Pi$ denote the Poisson tensor (bi-vector) on $M$ defined, as usual, by $\Pi(df,dg)=\{f,g\}$. A Hamiltonian vector field in the Poisson category is defined by the following equation $X_f:= \Pi(df, \cdot)$. The set of Hamiltonian vector fields on a Poisson manifold defines a distribution which is integrable in the sense that we may associate a (singular) foliation to it.
This foliation is called the \emph{symplectic foliation} of the Poisson manifold $(M,\Pi)$. The dimension of the symplectic leaf through a point coincides with the rank of the bivector field $\Pi$ at the point. The functions defining this foliation (if they exist) are called Casimirs. A manifold with constant rank is called a {\it regular} Poisson manifold.

The following theorem by Weinstein \cite{weinstein} gives the local picture of this foliation. A Poisson manifold is locally a direct product of a symplectic manifold, endowed with its Darboux symplectic form, together with a transverse Poisson structure whose Poisson vector field vanishes at the point.

\begin{Th}[Weinstein splitting theorem]
Let $(M^n,\Pi)$ be a Poisson manifold of dimension $n$ and let the rank of $\Pi$ be $2k$ at a point $p\in M$. Then, in a  neighbourhood  of $p$, there exists a coordinate system
$$
(x_1,y_1,\ldots, x_{k},y_{k}, z_1,\ldots, z_{n-2k})
$$
such that the Poisson structure can be written as
\begin{equation}
\Pi = \sum_{i=1}^k \frac{\partial}{\partial x_i}\wedge
\frac{\partial}{\partial y_i} + \sum_{i,j=1}^{n-2k}
f_{ij}(z)\frac{\partial}{\partial z_i}\wedge
\frac{\partial}{\partial z_j} ,
\end{equation}
where $f_{ij}$ are functions vanishing at the point and depending only on the transversal variables $(z_1,\ldots,z_{n-2k})$.
\end{Th}

The notion of integrable systems on a Poisson manifold (see also Section \ref{bi}) stands for a collection of integrable systems on this symplectic foliation. Namely,
  let $(M,\Pi)$ be a Poisson manifold of (maximal) rank $2r$ and of dimension $n$. An $s$-tuplet of Poisson-commuting functions
  $F=(f_1,\dots,f_s)$ on $M$ is said to define a \emph{Liouville integrable system} on $(M,\Pi)$ if $f_1,\dots,f_s$ are independent (i.e., their differentials are independent on a dense open set) and $r+s=n$.

The map $F:M\to\R^s$ is often called the \emph{moment map} of $(M,\Pi,F)$.

For integrable systems on Poisson manifolds it is possible to obtain an action-angle theorem in a  neighbourhood  of a regular torus, as it was proven in \cite{LMV}, showing that, up to a diffeomorphism, the set of integrals $F$ is indeed the moment map of a toric action:

\begin{Th}[Laurent-Gengoux, Miranda, Vanhaecke]
  Let $(M,\Pi)$ be a Poisson manifold of dimension $n$ of maximal rank $2r$. Suppose that $F=(f_1,\dots,f_s)$
  is an integrable system on $(M,\Pi)$ and $m\in M$ is a point such that
\begin{enumerate}
  \item[(1)] $df_1\wedge\dots\wedge df_s\neq0$;
  \item[(2)] The rank of $\Pi$ at $m$ is $2r$;
  \item[(3)] The integral manifold $L_m$ of $X_{f_1},\dots,X_{f_s}$, passing through $m$, is compact.
\end{enumerate}
  Then there exists ${\mathbb R} $-valued smooth functions $(\sigma_1,\dots, \sigma_{s})$ and $ {\mathbb R}/{\mathbb Z}$-valued smooth
  functions $({\theta_1},\dots,{\theta_r})$, defined in a  neighbourhood  $U$ of $L_m$, such that
  \begin{enumerate}

  \item The manifold $L_m$ is a torus $ T^r$;

    \item The functions $(\theta_1,\dots,\theta_r,\sigma_1,\dots,\sigma_{s})$ define an isomorphism
      $U\simeq T^r\times B^{s}$;
    \item The Poisson structure can be written in terms of these coordinates as
      $$
        \Pi=\sum_{i=1}^r\frac{\partial}{\partial \theta_i}\wedge\frac{\partial}{\partial \sigma_i},
      $$
      in particular, the functions $\sigma_{r+1},\dots,\sigma_{s}$ are Casimirs of $\Pi$ (restricted to $U$);
    \item The leaves of the surjective submersion $F=(f_1,\dots,f_{s})$ are given by the projection onto the
          second component $ T^r\times B^{s}$, in particular, the functions $\sigma_1,\dots,\sigma_{s}$ depend on the
          functions $f_1,\dots,f_{s}$ only.
  \end{enumerate}

\end{Th}

The functions $\theta_1,\dots,\theta_r$ are called {\it angle coordinates}, the functions $\sigma_1$,
$\dots,\sigma_r$ are called {\it action coordinates}, and the rest of the functions
$\sigma_{r+1},\dots,\sigma_{s}$ are called {\it transverse coordinates}.

Observe, in particular, that this theorem provides not only the existence of Liouville tori and action-angle coordinates in a neighbourhood of an invariant submanifold of a Poisson manifold, but also it ensures that, in a neighborhhod of an integral manifold $L_m$, there is a set of Casimirs. Difficulties arise when we try to extend this theorem to points where the Poisson structure is not regular. For particular classes of Poisson manifolds where some transversality conditions are met (such as  $b$-Poisson manifolds \cite{guimipi2, gmw}), it has been possible to extend this action-angle coordinates scheme (see \cite{km, kms} and also \cite{gmps}).

The first open problem concerns this extension:

\begin{problem}
Extend the action-angle theorem in a neighbourhood of a Liouville torus at singular points of the Poisson structure.
\end{problem}

The first adversity in doing so is that, contrary to the expectations, in general, an integrable system is not a product of two integrable systems (one on the symplectic leaf and the other one in the transverse part). So there is no simultaneous splitting theorem for the Poisson structure and the integrable system. Counterexamples to this decomposability can be found in \cite{camilleeva}. When this happens, we say the system is splittable.  So the problem above can be refined as follows:

\begin{problem}
Extend the action-angle theorem in a neighbourhood of a Liouville torus at singular points of the Poisson structure for splittable integrable systems.
\end{problem}

We suspect that either stability of the transversal part of the Poisson manifold or the existence of associated toric actions will play a role in solving this problem for special classes of Poisson manifolds.

This problem can be extended to its global version. Can we find action-angle coordinates globally? This problem was considered by Duistermaat in the symplectic case \cite{duistermaat}.  One of the main ingredients in his construction was the affine structure defined  by the action variables on the image of the moment map. In the regular Poisson setting, the image of the moment map is foliated by affine manifolds. Understanding this foliation is a key point to understanding global obstructions. The following problem has been partially considered in \cite{camilleruipol}.

\begin{problem}
Determine obstructions to the global existence of action-angle coordinates on regular Poisson manifolds.
\end{problem}

Realisation problem is also important:

\begin{problem}Which foliations with affine leaves can be described as the image of the moment map?
\end{problem}

Notice that a similar problem is still open even in the symplectic setting:  which closed manifolds $B^n$ may serve as  bases of locally trivial Lagrangian fibrations $M^{2n} \overset{T^n}{\longrightarrow} B^n$, where $M^{2n}$ is symplectic?

Extensions of the affine structure to the singular setting are also desirable. One instance when the affine structure on the leaves extends to an affine structure on the image of the moment map for non-regular Poisson manifolds is the case of $b$-Poisson manifolds \cite{gmps}. For these manifolds this affine structure is well-understood and classified via the Delzant theorem \cite{gmps2}. Such a phenomenon is also observed for other Poisson structures with open dense leaves on $b^m$-Poisson manifolds \cite{gmw3}.

\begin{problem} Consider a Poisson manifold $M$ of even dimension such that it is symplectic on a dense set $U\subset M$.  Assume that $M$ is  endowed with a toric action which is Hamiltonian on $U$. Is there an analog of Delzant theorem for the image of the moment map?
\end{problem}

Furthermore, in view of \cite{mirandascott} where the Moser homotopy method is presented also for foliated manifolds, can Delzant theorem be generalised for regular Poisson manifolds?

We conclude this section with a problem related to the symplectic case.  However, its solution might require some Poisson techniques.
There are two obstructions to globalizing action-angle coordinates, as
described by Duistermaat~\cite{duistermaat}: the Hamiltonian monodromy
and the Chern (or Duistermaat-Chern) class. Contrary to the former,
the latter has been very little investigated. There are, clearly, at
least two reasons for this:
\begin{itemize}
\item The Chern class always vanishes for two degrees of freedom systems;
\item There are no known physically meaningful examples with
  non-trivial Chern class in three degrees of freedom.
\end{itemize}
Although it is easy to manually construct any Chern class, see~\cite{bates-chern}, having a natural example from classical
mechanics would provide better motivation for further study (see also Problem~\ref{sec:define-and-detect}).

\begin{problem}[S. V\~u Ng{\d o}c] \label{sec:find-non-trivial}
Find natural examples of integrable systems in classical mechanics with non-trivial Duistermaat--Chern classes.
\end{problem}


\section{Integrability and Quantisation}

\subsection{Quantum Integrable Systems (communicated by S.V\~u Ng{\d o}c) }
\label{sec:quant-integr-syst}

The following open problems concern the quantum aspect of Liouville
integrable systems. From the quantum mechanical perspective, the most
natural notion of finite dimensional quantum integrable systems is
obtained by replacing classical observables (smooth real functions on
a symplectic manifold) by quantum observables (self-adjoint operators
on a Hilbert space), and, due to the correspondence principle of
Dirac, Poisson brackets by commutators.

Therefore, a quantum integrable system on a ``quantum Hilbert space
corresponding to $M^{2n}$'' should be the data of $n$
\emph{semiclassical} operators $P_1,\dots,P_n$ such that
\begin{enumerate}
\item \label{item:self-adjoint }$P_j^* = P_j$
\item \label{item:commute} $[P_i,P_j]=0$
\item \label{item:independent} the $P_j$'s are ``independent'' (see
  below, and Problem~\ref{sec:spectr-indep}).
\end{enumerate}

The word \emph{semiclassical} is very important and means that such
operators $P_j$ should have a ``classical limit'', typically in the
regime when Planck's constant $\hbar$ tends to zero. Without this
semiclassical limit, the number $n$, which relates to the dimension
$2n$ of the classical phase space, would be irrelevant.

There are several instances of quantum settings when the semiclassical
limit can be made precise. The two prominent theories are
pseudo-differential operators (when $M$ is a cotangent bundle) and
Berezin-Toeplitz operators (when $M$ is a prequantisable symplectic
manifold). In both cases, operators $P$ have a \emph{principal symbol}
$p=\sigma(P)$, which is a smooth function on $M$, and Dirac's
correspondence principle holds modulo $O(\hbar^2)$, which implies that
$[P_j,P_k]=0 \implies \{p_j,p_k\}=0$. We can now explain the last item
above: we say that the $P_j$'s are independent when the principal
symbols $p_j$'s have (almost everywhere) independent differentials.

The following very natural and classical question that one can ask immediately after
defining quantum integrable systems is still open and probably quite
difficult.

\begin{problem}
Assume that a classical integrable system $(p_1,\dots,p_n)$
is given on $M$. Does there exist a quantum integrable system
$(P_1,\dots,P_n)$ such that $p_j=\sigma(P_j)$?
\end{problem}

See Charles {\it et al} \cite{ChPeVN13} for a solution in the compact toric case.
More specific versions of this general question are discussed in Section \ref{quantum}.

The difficulty is to construct operators $P_j$ with the \emph{exact}
commutation property $[P_j,P_k]=0$. Using microlocal analysis and, in
particular, Fourier integral operators~\cite{FIO2}, one should be able
to obtain $[P_j,P_k]=O(\hbar^\infty)$. In many classical examples, an
exact quantisation is known, see for instance~\cite{toth-quadratic,
  heckman-quantum}. Symmetry considerations should help greatly. This
is a typical local-to-global obstruction problem, and one expects a
corresponding cohomology complex. In the case of polynomial symbols in
$\mathbb C^n$, the complex was constructed
in~\cite{garay-vanstraten-integrability}.

\begin{problem}[S. V\~u Ng{\d o}c]\label{sec:spectr-indep}
Given a set of semiclassical operators that verify
conditions~\ref{item:self-adjoint } and \ref{item:commute} above, can
one detect the independence axiom~\ref{item:independent} in a purely
spectral way? More precisely, can one tell from asymptotics of joint
eigenvalues or eigenfunctions of $(P_1,\dots,P_n)$ that the principal
symbols are almost everywhere independent?
\end{problem}

A great achievement of microlocal analysis of pseudo-differential
operators is the formulation of   \emph{complete Bohr-Sommerfeld rules
  to any order}. Given a quantum integrable system, these rules are, in
fact, some kind of integrality conditions (or cohomological conditions),
depending on a parameter $E\in\R^n$, that are satisfied \emph{if and
  only if} $E$ is a joint eigenvalue of the system, modulo an
arbitrarily small error of size $O(\hbar^\infty)$.

Historically (both in physics and mathematics), these rules have been
written for regular Lagrangian tori, {i.e.}, $E$ should stay in a
small ball of regular values of the classical moment map. More
recently, non-trivial extensions to some singular situations have been
found, see~\cite{child-book, san-panoramas, san-techniques}.

In order to treat even more physical examples, and to appeal to the
geometry community, it is important to extend these rules to
Berezin-Toeplitz quantisation on compact K\"ahler manifolds. This was
done recently for regular tori~\cite{Ch2003a}, and for one degree
of freedom systems with Morse singularities,
see~\cite{lefloch-elliptic,lefloch-hyperbolic}. It should be possible
to treat the focus-focus case as well, along the lines
of~\cite{san-focus}.

\begin{problem}[S. V\~u Ng{\d o}c] \label{sec:write-bohr-somm}
Write Bohr-Sommerfeld rules for focus-focus singularities in Berezin-Toeplitz quantisation.
\end{problem}

A similar problem appears in the context of geometric quantisation (see Problem \ref{EvaFocus} below).

In relation to Problem~\ref{sec:find-non-trivial}, the following
natural question is still open: what is the quantum analogue of the
Chern (or Duistermaat-Chern) class? Contrary to Quantum Monodromy,
which was defined purely in terms of the joint
spectrum~\cite{cushman-duist, san-mono}, it is quite possible that the
Chern class manifests itself as a phase on eigenfunctions, and thus
could be invisible in the spectrum. This study probably involves
``shift operators'' that are transversal to the Lagrangian fibration
and allow to ``jump'' from one eigenspace to another one.

\begin{problem}[S. V\~u Ng{\d o}c]\label{sec:define-and-detect}
Define (and detect) the quantum Chern class.
\end{problem}

\subsection{Inverse problem for quantum integrable systems (communicated by \'{A}. Pelayo and S. V\~u Ng{\d o}c) }

The notion of quantum integrable system as a maximal set of commuting quantum observables
is fairly well established now. It is a very old notion, going back to  Bohr, Sommerfeld and
Einstein~\cite{E1917} in the early days of quantum mechanics. Nonetheless the theory didn't
blossom at that time because of the lack of technical tools needed to solve interesting problems
about quantum integrable systems.  For recent surveys discussing the spectral geometry of
integrable systems, see~\cite{DCDS12, Pe13}.

The most elementary results about the symplectic structure of classical integrable systems, such
as the existence of action-angle variable, did not fit well in Schr{\"o}dinger's quantum setting at that
time because they involve the analysis of differential and pseudodifferential operators, only
developed since the 1960s. This is a prominent part of analysis and PDE (the most mathematical side of
quantum mechanics), that goes by the name of ``microlocal analysis".

The microlocal analysis of action-angle variables starts with the pioneering works of
Duistermaat~\cite{Du1974} and Colin de Verdi{\`e}re \cite{CdV,CdV2}, and was followed by
many other authors (see~\cite{ChPeVN13} for further references). For compact symplectic manifolds,
the study of quantum action-angle variables is very recent, it goes back to Charles \cite{Ch2003a},
using the theory of Toeplitz operators.   Toeplitz operators give rise to a semiclassical algebra of
operators with symbolic calculus and microlocalization properties,  isomorphic (at a microlocal level)
to the algebra of pseudodifferential operators (see Boutet de Monvel\--Guillemin~\cite{BoGu1981}).

The classification of semitoric systems
in~\cite{Pelayo2, Pelayo3} (see Section 1)  goes a long way to answering
the fascinating question ``can one hear the shape of a semitoric
system''? More precisely, the idea is to recover the classical
semitoric invariants from the joint spectrum of a quantum integrable
system, as much as possible.

\begin{con}[Pelayo-V\~u Ng\d oc 2011, Conjecture 9.1 in \cite{Pelayo}]  \label{PeVNconjecture}
A semitoric system $F=(J, H)$ is determined, up to symplectic equivalence, by its semiclassical joint
spectrum (i.e., the collection of points $(\lambda,\nu) \in \mathbb{R}^2$, where $\lambda$ is a
eigenvalue of $\hat{J}_{\hbar}$ and $\nu$ is an eigenvalue of $\hat{H}_{\hbar}$, restricted to the
$\lambda$\--eigenspace of $\hat{J}$, as $\hbar \to 0$). Moreover, from any such semiclassical
spectrum, one can explicitly construct the associated semitoric system.
\end{con}

This conjecture is now a theorem for the special cases of quantum
toric integrable systems on compact manifolds~\cite{ChPeVN13}, and for semitoric integrable
systems of Jaynes\--Cummings type~\cite{LFPeVN2016}. The problem is that it is unclear whether
one of the most subtle (perhaps the most subtle) invariants of integrable systems, known as the twisting
index, can be detected in the semiclassical spectral asymptotics. If this is the case (which is very likely),
then the  Pelayo-V\~u Ng\d oc inverse spectral conjecture
(\cite{Pelayo}) above will hold.

An important part of the problem is to recover the Taylor series
invariant that classifies semi-global  neighbourhood s of focus-focus
fibers~\cite{San}. A positive solution was given
in~\cite{PeVN12-CMP, LFPeVN2016}, but the
procedure was not constructive: one had to first let $\hbar\to 0$ for a
regular value $c$ near the focus-focus value $c_0$, and then take the
limit $c\to c_0$, which doesn't help computing the Taylor series in an
explicit way from the spectrum.

\begin{problem}[S. V\~u Ng{\d o}c]\label{semitoric_drum}
Compute the Taylor series invariant of semitoric systems directly from the spectrum ``at'' the focus-focus
critical value.
\end{problem}

One should probably use the corresponding singular
Bohr-Sommerfeld rule,  see Problem~\ref{sec:write-bohr-somm}. In fact,
the following intermediate question is very natural:

\begin{problem}[S. V\~u Ng{\d o}c] Is the singular Bohr-Sommerfeld formal power series in $\hbar$ a spectral invariant?
\end{problem}

There are various later refinements and extensions of
Conjecture \ref{PeVNconjecture}~\cite{DCDS12, PePoVN, SeVN}, following other inverse
spectral questions in semiclassical analysis (see for
instance~\cite{zelditch-inverse-II, hezari,
  iantchenko-sjostrand-zworski})
which concern integrable systems (or even collections of
commuting operators) more general than semitoric, and even in higher dimensions. All of them essentially
make the same general claim: ``from the semiclassical  joint spectrum of a quantum integrable system
one can detect the principal symbols of the system".  A more concrete problem is:

\begin{problem}[\'{A}. Pelayo]
What information about the principal symbols $f_1,\ldots, f_n$ of
a quantum integrable system $T_{1,\hbar},\ldots, T_{n,\hbar}$ can be detected from their semiclasical
joint spectrum? For instance, suppose that we know that a certain object -- say, an integer $z$, or a
matrix $A$ -- is a symplectic invariant of $(M, \omega, f_1,\ldots, f_n)$. Can we compute this object
from the semiclassical joint spectrum?
\end{problem}

An additional interesting problem, which is in fact closely related to Question \ref{prob:1.4}  is as follows:

\begin{problem}[\'{A}. Pelayo] Can one detect from the joint spectrum of a quantum integrable system $T_{1,\hbar},\ldots, T_{n,\hbar}$ that a singularity is degenerate?
\end{problem}

The theory of semitoric systems provides a strong motivation for understanding  the connections between
integrable systems and Hamiltonian torus actions. In view of this,  we conclude by discussing a possible
quantum approach to counting fixed points of symplectic $S^1$\--actions.

In view of Tolman's recent work constructing
a non-Hamiltonian symplectic $S^1$\--action on a compact manifold
with isolated fixed points, the questions
 how many fixed points symplectic actions have (both upper and lower bounds) and  what the relation  between having fixed points and being Hamiltonian is,
have attracted much interest~\cite{ToWe, Go, PeTo11, GoPeSa}.  See also Pelayo \cite{AP} for a discussion and other references.   As far as we know, the ``quantised"
version of this problem is not well  understood, although some steps have been taken in this direction~\cite{LFPesemiclassicallimits}.
The object to quantise can no longer be the momentum map, since  there is no momentum map --- however, McDuff pointed
out in~\cite{Mc88} that any symplectic $S^1$\--action admits an  $S^1$\--valued momentum map $\mu \colon M \to S^1$.
The explicit construction  of this map appears in ~\cite{PeRa12},  and it can be quantised~\cite[Section 7]{LFPesemiclassicallimits}
(here the operator will be unitary, rather than self\--adjoint).

\begin{problem}[\'{A}. Pelayo]
 Can one make progress in counting the number of fixed points by studying the spectrum of the quantisation of $\mu \colon M \to S^1$?
\end{problem}

This type of problem should be very relevant to the study of equilibrium points of integrable systems but, for the moment,
it is only a toy version to start with. The general case would involve quantisation and a study of the fixed point properties
of an object (momentum map) of the form $M \to (S^1)^k \times \mathbb{R}^n$, and this will likely be quite challenging.

\subsection{Quantum integrability for polynomial in momenta integrals}  \label{quantum}

The approach to quantisation of natural Hamiltonian systems discussed in this section is, perhaps, the most traditional.
Instead of  a Hamiltonian of the form (\ref{ham}), i.e., $H = K+U$ on a Riemannian manifold $M$, one  considers the second order differential operator acting on $C^2$-differentiable functions on $M$ and given by
$$
\hat H:= \Delta + U,
$$
where $\Delta= g^{ij} \nabla_i \nabla_j$ is the Beltrami-Laplace operator corresponding to $g$. Next, by a quantum  integral we  understand another differential operator $\hat F$  that  commutes with $H$:
\begin{equation} \label{commute}
[\hat H, \hat F]= \hat H \circ \hat F - \hat F \circ \hat H=0 .
\end{equation}

Sometimes one requires that the operator $\hat F$ be self-adjoint.
It is known and is easy to see that, for any operators $ \hat H $ and $ \hat F $,  the symbol of  $[\hat H, \hat F]$ is just the Poisson bracket of the symbols of $\hat H$ and $ \hat F$.  We say that a polynomial in momenta
integral $F$ for (\ref{ham}) is {\it quantisable} if there exists a quantum integral $\hat F$ for $\hat H$ such that its symbol is $F$.

 \begin{problem}  \label{MatveevQuant}
 What are necessary and/or sufficient conditions on a metric $g$ such that every  polynomial integral of its geodesic flow is quantisable?
\end{problem}

Partial answers are known in the case of integrals of small degrees and in dimension 2.

A  more specific version of Problem \ref{MatveevQuant} is related to  quantisation of the argument shift method (see Section \ref{argshift}).  In this setting, we impose an additional condition on classical and quantum integrals.  Namely, we assume that $M=G$ is a Lie group, and the polynomial commuting integrals are all left invariant.   Is it possible to find the corresponding commuting quantum integrals, i.e., differential operators, that are still left invariant?   More specifically, this question is related to a special family (algebra) of left invariant polynomial integrals on $T^*G$ obtained by the argument shift method  \cite{MF78-Izv}, and it can be asked in a purely algebraic terms.

Let $\mathcal P(\Fg)$ be the algebra of polynomials on $\goth g^*$ endowed with the standard Lie-Poisson bracket \eqref{LiePoisson}, and $a\in \goth g^*$ be a regular element.  For a local $\mathrm{Ad}^*$-invariant function $f$ (not necessarily polynomial), consider its Taylor expansion at the point $a\in\goth g^*$:
$$
f(a + \lambda x) = f(a) + \lambda f_1(x) + \lambda^2 f_2(x) + \dots + \lambda^k f_k(x) + \dots
$$
with $f_k(x) \in \mathcal P(\Fg)$ being a homogeneous polynomial of degree $k$.  The subalgebra $\mathcal F_a \subset \mathcal P(\goth g)$, generated by  $f_k$'s  (where $k\in \mathbb N$ and $f$ runs over the set of local $\mathrm{Ad}^*$-invariant functions), is commutative with respect to the Lie-Poisson bracket \eqref{LiePoisson} and is called a {\it Mischenko-Fomenko subalgebra} of $\mathcal P(\goth g)$.  Each polynomial $g\in \mathcal F_a$ can be naturally lifted to the cotangent bundle $T^*G$ and, as a result, we obtain a commutative subalgebra of left-invariant functions, polynomial in momenta on $T^*G$.  Quantising this algebra in the class of left-invariant differential operators is equivalent to ``lifting''  $\mathcal F_a$ to the universal enveloping algebra $U(\mathfrak g)$, i.e., constructing a commutative subalgebra $\widehat{\mathcal F}_a\subset U(\goth g)$ such that for any element $\hat p \in \widehat{\mathcal F}_a$ of degree $k$, its ``principal symbol''
$p \in \mathcal P(\goth g)$ (being a homogeneous polynomial of degree $k$) lies in the algebra  $\mathcal F_a$, and $\mathcal F_a$ is generated by such elements.

 A method for constructing the quantum Mischenko-Fomenko algebra  $\widehat{\mathcal F}_a$ for a semisimple
Lie algebra $\mathfrak g$ was suggested by L.G.\,Rybnikov~\cite{Rybnikov}. However, it is not clear whether  a similar construction works for an arbitrary Lie algebra $\mathfrak g$.

\begin{problem}
Quantise the Mischenko-Fomenko algebra   $\mathcal F_a \subset \mathcal P(\goth g)$ for an arbitrary finite-dimensional Lie algebra $\mathfrak g$ {\/\rm(}or find an obstruction to  quantisation{\/\rm)}.
\end{problem}

\subsection{Integrable systems and geometric quantisation}

As explained above,  the idea of quantisation is to associate a representation space with a classical system  in a functorial way, so that the observables (smooth
functions) become operators on a Hilbert space, and the classical
Poisson bracket becomes the commutator of operators. Here we would like to stress the role of Lagrangian foliations in Geometric Quantisation.

Consider a symplectic manifold $(M^{2n},\omega)$ with an integral homology class $[\omega]$. Because of integrality of  $[\omega]$ (\cite{woodhouse},
\cite{kirillov}), there exists a complex line bundle $\mathbb{L}$
with a connection $\nabla$ over $M$ such that $curv(\nabla)=\omega$.  The pair $(\mathbb{L},\nabla)$ is called {\emph{a prequantum line
bundle}} of $(M^{2n},\omega)$.  Having fixed these data, as a first candidate for the
representation space, we consider the set of sections  $s$  which are flat, meaning that $\nabla_X s=0$ in some privileged directions,
 given by polarisation (integrable Lagrangian distribution of the complexified tangent bundle).  Any Lagrangian foliation of $M^{2n}$ defines a polarisation,  and thus integrable systems can be seen as real polarisations,  and their tangent directions can be used to solve the flatness equation $\nabla_X s=0$.

 It is convenient to extend the notion of polarisation to take singularities into account. Thus a real polarisation $\mathcal P$ is a
foliation  whose generic leaves are Lagrangian submanifolds.

A leaf of a polarisation is called \emph{Bohr-Sommerfeld} if the flat sections are defined along the leaf. In the classical case, when $M^{2n}$ is a cotangent bundle with the leaves being the cotangent spaces, this condition is fulfilled for all of them. However, if the leaves of a polarisation are compact, then this condition becomes nontrivial.

As we will see next, the action-angle coordinates of integrable systems are useful for localisation of Bohr-Sommerfeld leaves, and a simple model of quantisation would help to count them. This coincides with Kostant's viewpoint on geometric quantisation (see \cite{kostant}).

As the simplest example,  consider the surface $M^2=S^1\times \mathbb
R$ endowed with the symplectic form $\omega=dt\wedge d\theta$. Take, as a prequantum line bundle $\mathbb L$, the trivial
bundle with connection $1$-form $\Theta= td\theta$, and as a distribution $\mathcal
P=\langle \frac{\partial}{\partial \theta}\rangle$.
The  flat sections satisfy
$\nabla_X\sigma=X(\sigma)-i\langle \theta,X \rangle\sigma=0$.

 Thus
$\sigma(t,\theta)=a(t) e^{it\theta}$, and Bohr-Sommerfeld leaves are
given by the condition $t=2\pi k, k\in \mathbb Z$, i.e., by the integral action coordinates, associated to the polarisation with the moment map $t$.

 This characterisation of
Bohr-Sommerfeld leaves generalises  for regular fibrations:
\begin{Th}[Guillemin-Sternberg, \cite{guilleminandsternberg}] Assume the polarisation is a
regular Lagrangian fibration with compact leaves over a simply connected base
$B$. Then the Bohr-Sommerfeld set is discrete;  and if the
zero-fiber  $\{(f_1(p),\dots, f_n(p))=0\in \R^n\} \subset M^{2n}$ is a  Bohr-Sommerfeld leaf, then the Bohr-Sommerfeld set is
given by $BS=\{p\in M, (f_1(p),\dots, f_n(p))\in\mathbb Z^n\}$,
where $f_1,\dots,f_n$ are global action coordinates on $B$.
\end{Th}

Let us now relate Bohr-Sommerfeld leaves with the geometric quantisation scheme.
 Denote by $\mathcal{J}$  the sheaf
of (local) flat sections along the polarisation. Following Kostant \cite{kostant},  we define the
quantisation as the sheaf cohomology  $\mathcal{Q}(M)=\bigoplus_{k\geq
0}H^k(M,\mathcal{J})$. The following theorem by Sniatycki \cite{sniatpaper} provides a bridge between  this sheaf cohomology approach of Kostant and counting Bohr-Sommerfeld leaves:

\begin{Th}[Sniatycki] If the leaf space $B^n$ is Haussdorf and the natural projection $\pi:M^{2n}\to B^n$ is
a fibration with compact fibers, then all the cohomology groups $H^k(M,\mathcal{J})$
vanish,  except for $k=\frac{1}{2}\dim M$.
Thus, $\mathcal{Q}(M^{2n})=H^n(M^{2n},\mathcal{J})$, and the
dimension of the vector space $H^n(M^{2n},\mathcal{J})$ equals the number of
Bohr-Sommerfeld leaves.
\end{Th}

This scheme works well for toric manifolds, where the base $B$  may be identified with the image of the
moment map associated with the torus action, i.e., the interior of the Delzant polytope  and, more generally, for Lagrangian fibrations with elliptic singularities \cite{hamilton}.    The next step would be to generalise this scheme to the case of integrable systems with arbitrary non-degenerate singularities.
Unfortunately, hyperbolic and focus-focus singularities \cite{hamiltonmiranda}, \cite{mirandapresassolha} do not behave well under this scheme, and the dimension of $H^n(M^{2n},\mathcal{J})$ may become infinite  even if the number of Bohr-Sommerfeld leaves is finite.

\begin{problem}[Miranda-Presas-Solha \cite{mirandapresassolha}]\label{EvaFocus}
  Modify this scheme to get finite dimensional representation spaces for focus-focus and hyperbolic singularities that still capture the invariants of the corresponding integrable systems and their singularities.
\end{problem}


\end{document}